\documentclass[12pt]{article}

\usepackage{latexsym}
\newtheorem{theoreme}{{\bf Th\'eor\`eme}}[section]

\newtheorem{corollaire principal}[principal]{\bf Corollaire}
\newtheorem{proposition}[theoreme]{{\bf Proposition}}
\newtheorem{lemme}[theoreme]{{\bf Lemme}}
\newtheorem{corollaire}[theoreme]{{\bf Corollaire}}
\newtheorem{definition}[theoreme]{{\bf D\'efinition}}

\newenvironment{demonstration}{\noindent{\bf D\'emonstration
}}{\nolinebreak $\Box $\hspace{-2.15mm}\raisebox{1.25mm}{.} \medskip}

\newenvironment{demonstration du lemme}{\noindent{\bf D\'emonstration du lemme
}}{\nolinebreak $\Box $\hspace{-2.15mm}\raisebox{1.25mm}{.} \medskip}

\def\CC{{\bf C}}

\def\NN{{\bf N}}

\begin{document}

\title{{\Large \bf  Homog\'en\'eit\'e locale pour les m\'etriques riemanniennes holomorphes en dimension $3$}}

\author{{\normalsize
{\bf Sorin DUMITRESCU}}}
\date{F\'evrier   2006}

\maketitle

\vspace{0.3cm}
{\normalsize{{\bf Homog\'en\'eit\'e locale pour les m\'etriques riemanniennes holomorphes en dimension $3$}}}

\
\vspace{0.1cm}

\noindent{\bf{ R\'esum\'e.}} 
Nous d\'emontrons que sur les vari\'et\'es complexes compactes de dimension $3$ les m\'etriques riemanniennes holomorphes sont n\'ecessairement localement homog\'en\'es (i.e. le
pseudo-groupe des isom\'etries locales agit transitivement sur la vari\'et\'e.)

\section{Introduction}
 
 Dans le contexte de la g\'eom\'etrie complexe (holomorphe),  une   {\it m\'etrique riemannienne holomorphe} est l'analogue du concept r\'eel de m\'etrique pseudo-riemannienne.
 Formellement une m\'etrique riemannienne holomorphe sur une vari\'et\'e complexe $M$ de dimension $n$ est une section holomorphe $q$ du fibr\'e $S^2(T^{*}M)$  des formes quadratiques complexes sur l'espace tangent holomorphe \`a $M$ telle  que, en tout point $m$ de $M$, la forme quadratique complexe $q(m)$ est non d\'eg\'en\'er\'ee (de rang maximal, \'egal \`a la dimension complexe de $M$).
 
 L'exemple type d'une m\'etrique riemannienne
 holomorphe est la m\'etrique {\it plate}  $q=dz^1 + dz^2 + \ldots + dz^n$ sur $\CC ^n$.  A l'instar du cadre pseudo-riemannien (voir, par exemple,~\cite{Wo}), il n'y a pas de difficult\'e \`a introduire le tenseur de courbure 
  d'une m\'etrique riemannienne holomorphe  comme obstruction infinit\'esimale (\`a l'ordre $2$) \`a la platitude de la m\'etrique.  Il existe une unique connexion lin\'eaire  holomorphe  
  (de Levi-Civita) sur $M$ compatible avec la m\'etrique riemannienne holomorphe et cette connexion permet de d\'efinir des courbes g\'eod\'esiques qui g\'en\'eralisent les droites affines (parcourues \`a vitesse constante) du cas plat (voir, par exemple,~\cite{Le} pour une \'etude des g\'eod\'esiques holomorphes).
  
  Nous allons nous int\'eresser dans cet article aux m\'etriques riemanniennes holomorphes  sur les vari\'et\'es complexes {\it compactes.} L'existence d'une telle m\'etrique sur
  une vari\'et\'e complexe compacte impose des conditions tr\`es restrictives \`a la vari\'et\'e. Une premi\`ere obstruction \'evidente est la premi\`ere classe de Chern. En effet,
  dans le langage des $G$-structures d\^u \`a C. Ehresmann (la terminologie \'etant due \`a S. Chern) une m\'etrique riemannienne holomorphe est une $O(n, \CC)$-structure sur le fibr\'e tangent holomorphe \`a $M$. Ceci implique l'existence d'un rev\^etement double non ramifi\'e de $M$ sur lequel  le groupe structural du fibr\'e des rep\`eres (qui est, en g\'en\'eral, un $GL(n, \CC)$-fibr\'e) se r\'eduit au groupe
  $SO(n, \CC)$ qui pr\'eserve un volume holomorphe. Le fibr\'e canonique (des formes volumes holomorphes) de ce rev\^etement double de $M$  est donc trivial, ce qui assure
  l'annulation de la premi\`ere classe de Chern de $M$. On peut aussi remarquer que la pr\'esence d'une m\'etrique riemannienne holomorphe fixe un isomorphisme
  entre le fibr\'e tangent holomorphe $TM$ et son dual $T^{*}M$. En particulier, le fibr\'e canonique de $M$ est isomorphe \`a son dual ce qui implique l'annulation de la premi\`ere
  classe de Chern de $M$.
  
  Un premier exemple de vari\'et\'es compactes admettant des m\'etriques riemanniennes holomorphes est donn\'e par les tores complexes.  Pour s'en assurer il suffit de constater 
  que la m\'etrique plate $q=dz^1 + dz^2 + \ldots + dz^n$ est invariante par translation et descend donc sur tout quotient de $\CC ^n$ par un r\'eseau de translations.
  
  Pour ce qui est du cadre {\it k\"ahl\'erien} un r\'esultat de~\cite{IKO} montre que l'exemple pr\'ec\'edent est le seul (\`a rev\^etement fini pr\`es).  Plus pr\'ecis\'ement, il est d\'emontr\'e
  dans~\cite{IKO} que  parmi les vari\'et\'es k\"ahl\'eriennes compactes, seuls les tores complexes et leurs
 quotients finis admettent des connexions lin\'eaires 
  holomorphes (et donc
 des m\'etriques riemanniennes holomorphes). Dans ce cas, les m\'etriques riemanniennes holomorphes
 sont n\'ecessairement plates et donc, en particulier, {\it localement homog\`enes}, autrement dit le pseudo-groupe des biholomorphismes locaux qui pr\'eservent la m\'etrique riemannienne holomorphe (isom\'etries locales) agit transitivement sur $M$. 
  
 Un r\'esultat similaire est d\'emontr\'e dans~\cite{D2} pour les surfaces complexes compactes~: seuls les tores complexes et leurs
 quotients finis  admettent des m\'etriques riemanniennnes holomorphes et ces m\'etriques sont n\'ecessairement plates.
 
 Mentionnons \'egalement, dans le cadre k\"ahl\'erien,  les r\'esultats de sym\'etrie et classification
 obtenus pour des structures g\'eom\'etriques holomorphes plus g\'en\'erales  sur les vari\'et\'es de Calabi-Yau~\cite{D1}  et  sur les vari\'et\'es
 projectives unirationnelles~\cite{HM}.
 
 Nous nous int\'er\'essons dans cet article aux vari\'et\'es complexes compactes non (n\'ec\'essairement) k\"ahl\'eriennes de dimension complexe $3$ qui poss\`edent des m\'etriques
 riemanniennes holomorphes.
 
 Comme nous allons le voir dans la section 2, dans le cas  non K\"ahler, des exemples de vari\'et\'es complexes compactes poss\`edant des m\'etriques
 riemanniennes holomorphes sont les vari\'et\'es parall\'elisables, quotiens d'un groupe de Lie complexe par un r\'eseau co-compact.  Ces exemples  n'ont rien de
 surprenant car le fibr\'e tangent \'etant holomorphiquement trivial les m\'etriques riemanniennes holomorphes sur ces vari\'et\'es sont ''constantes'' par rapport au
 rep\`ere mobile invariant par translations (qui constitue  une trivialisation du fibr\'e tangent).
 
 En dimension $3$ des exemples in\'edits de vari\'et\'es complexes admettant des m\'etriques riemanniennes holomorphes ont \'et\'e construits dans~\cite{Gh} par E. Ghys.
 
 N\'eanmoins dans tous les exemples connus le rev\^etement universel de la vari\'et\'e est un groupe de Lie complexe et la pr\'eimage de la m\'etrique riemannienne holomorphe
 par le rev\^etement est invariante par les translations du groupe.

 Nous d\'emontrons ici le th\'eor\'eme suivant  qui g\'en\'eralise et pr\'ecise sous une forme optimale un premier r\'esultat  que nous avions obtenu dans~\cite{D2}.

\begin{theoreme}   \label{homogene} : Soit $M$ une vari\'et\'e complexe compacte connexe de dimension $3$ munie d'une m\'etrique riemannienne holomorphe. Alors toute structure g\'eom\'etrique holomorphe
de type affine sur M est n\'ecessairement localement homog\`ene. 

En particulier, la m\'etrique riemannienne holomorphe est localement homog\`ene.

\end{theoreme}

Il convient de situer ce r\'esultat dans le cadre de l'\'etude des {\it structures g\'eom\'etriques rigides} dans le sens de M. Gromov initi\'e dans~\cite{Gro} (voir \'egalement l'expos\'e de survey~\cite{DG} et~\cite{Ben}).
Sans donner de d\'efinition g\'en\'erale, pr\'ecisons que les m\'etriques riemanniennes holomorphes, ainsi que les m\'etriques pseudo-riemanniennes (dans le contexte r\'eel) sont des exemples type
de telles structures g\'eom\'etriques rigides. Une conjecture vague \'enonc\'ee par M. Gromov affirme que la pr\'esence d'un ''gros'' groupe d'isom\'etries pour une structures g\'eom\'etrique rigide sur une vari\'et\'e
compacte doit repr\'esenter une situation  suffisamment sym\'etrique (rare)  pour \^etre classifiable. En principe, il ne devrait y avoir que des exemples {\it alg\'ebriques} (construits \`a partir
de groupes de Lie) et leurs eventuels avatars (d\'eformations etc).  Les exemples abondent dans ce sens en g\'eom\'etrie riemannienne ou
pseudo-riemannienne r\'eelle (voir par exemple~\cite{Zeg}, pour une concr\'etisation de ce ph\'enom\`ene en g\'eom\'etrie lorentzienne de dimension $3$). Le  th\'eor\`eme~\ref{homogene}  vient conforter cette conjecture dans le cadre des vari\'et\'es complexes. Ici le caract\`ere holomorphe remplace l'hypoth\`ese
d'existence d'un groupe important d'isom\'etries et suffit pour engendrer des sym\'etries (isom\'etries) locales. 

Sans donner pour le moment la d\'efinition g\'en\'erale de structure g\'eom\'etrique (pour cela le lecteur devra se rapporter \`a la section $3$), pr\'ecisons  que dans le cas o\`u
la vari\'et\'e $M$ munie d'une  m\'etrique riemannienne holomorphe poss\`ede par ailleurs un (autre) tenseur  holomorphe $Y$ (par exemple, un champ de vecteurs) le
th\'eor\`eme~\ref{homogene} affirme que le pseudo-groupe des isom\'etries locales de la m\'etrique riemannienne holomorphe qui pr\'eservent en m\^eme temps le tenseur
$Y$ agit transitivement sur $M$. Autrement dit, m\^eme la structure g\'eom\'etrique ''totale''  (la plus riche possible) qui englobe \`a la fois la m\'etrique riemannienne holomorphe initiale et toute autre structure g\'eom\'etrique holomorphe globale de type affine  (voir la section $3$ pour la d\'efinition) situ\'ee sur $M$ est localement homog\`ene.

Nos m\'ethodes m\'elangent \`a la fois les arguments de g\'eom\'etrie diff\'erentielle rigide et des techniques qui viennent de la  g\'eom\'etrie analytique complexe. L'homog\'en\'eit\'e locale
d\'et\'ect\'ee par le th\'eor\`eme~\ref{homogene} conduit dans certains cas \`a des r\'esultats de classification~:

\begin{corollaire}  \label{corollaire} Soit $M$ une vari\'et\'e complexe compacte connexe de dimension $3$ munie d'une m\'etrique riemannienne holomorphe. Si $M$ poss\`ede une structure g\'eom\'etrique
holomorphe de type affine et de type g\'en\'eral, alors $M$ admet un rev\^etement fini non ramifi\'e qui est un quotient d'un groupe de Lie complexe connexe et simplement connexe de dimension $3$ par un r\'eseau
co-compact.

\end{corollaire}

Conjecturalement le r\'esultat du corollaire pr\'ec\'edent devrait rester valide en g\'en\'eral :
le rev\^etement universel d'une vari\'et\'e complexe compacte (de dimension 3) poss\`edant une
m\'etrique riemannienne holomorphe est un groupe de Lie complexe sur lequel la pr\'eimage
de la m\'etrique riemannienne holomorphe est invariante par translations. Le th\'eor\`eme
principal de cet article ouvre la voie vers ce type de r\'esultats de classification. Il s'agit  ensuite
de prouver des th\'eor\`emes dites {\it de compl\'etude} et d'arriver \`a int\'egrer les isom\'etries locales
obtenues par le th\'eor\`eme~\ref{homogene} en un groupe d'isom\'etries globales du rev\^etement universel qui agit transitivement. La question de la compl\'etude est ouverte m\^eme dans
le cas o\`u la m\'etrique riemannienne holomorphe est suppos\'ee plate (de courbure sectionnelle nulle) :
dans ce cas   il s'agit de prouver que le rev\^etement universel de $M$ est n\'ecessairement
biholomorphe \`a $\CC^3$.

Apr\`es cette courte introduction, la composition de cet article est la suivante.
La deuxi\`eme section  pr\'esente des constructions de vari\'et\'es complexes compactes poss\`edant des m\'etriques riemanniennes holomorphes. Dans la troisi\`eme
section nous rappelons bri\`evement les concepts  et les th\'eor\`emes n\'ecessaires \`a l'\'etude des structures g\'eom\'etriques rigides et de leurs isom\'etries locales. La section quatre 
est consacr\'ee \`a la d\'emonstration d'une forme faible du th\'eor\`eme principal :  on a homog\'en\'eit\'e locale sur un ouvert dense (en dehors d'un sous-ensemble analytique compact $S$ inclus dans $M$, eventuellement vide). Dans la derni\`ere section nous montrons comment on peut \'etendre cet ouvert dense \`a $M$ tout entier, ce qui ach\`eve la preuve  du th\'eor\`eme principal de l'article.

 Je remercie chaleureusement \'Etienne Ghys  pour m'avoir introduit au sujet  et pour son soutien constant. Je remercie \'egalement Ghani Zeghib pour la gentillesse avec laquelle il a toujours repondu \`a mes questions. Mes coll\`egues Nicolas Bergeron, Fran\c cois B\'eguin, Charles Frances et Pierre Pansu m'ont souvent \'ecout\'e parler de ce travail :  je les remercie pour leur patience et pour la pertinence de leurs suggestions.

  \section{Exemples}

  Les exemples les plus simples de vari\'et\'es complexes compactes qui admettent des m\'etriques riemanniennes holomorphes
  sont les vari\'et\'es dites parall\'elisables qui s'ob\-tienn\-ent comme quotient d'un groupe de Lie complexe par un r\'eseau
  co-compact.
  
  Un proc\'ed\'e simple pour construire de telles m\'etriques est le suivant : consid\'erons un groupe de Lie complexe $G$
  qui admet un r\'eseau co-compact $\Gamma$ (exemple : un groupe de Lie complexe semi-simple comme $SL(2,\CC)$).
  
   Nous identifions son alg\`ebre de Lie $\mathcal G$ \`a l'alg\`ebre de Lie des champs de vecteurs
     sur $G$ invariants par les translations \`a droite. Le quotient \`a droite par
    l'action de $\Gamma$ fournit une vari\'et\'e complexe compacte parall\'elisable $M=G/\Gamma$, 
    dont le fibr\'e tangent est holomorphiquement isomorphe \`a $M \times \mathcal G$. L'ensemble des m\'etriques holomorphes $q$
    sur $M$ s'identifie aux formes quadratiques non d\'eg\'en\'er\'ees
    sur $\mathcal G$. 
    
    Au cas o\`u cette forme quadratique est invariante par la repr\'esentation adjointe de $G$ dans son
    alg\`ebre de Lie (ce qui arrive, par exemple,  pour la forme de Killing d'un groupe de Lie
    semi-simple complexe : la forme de Killing, qui est toujours invariante par la repr\'esentation adjointe, est
    non d\'eg\'en\'er\'ee si $G$ est semi-simple) les translations \`a gauche sont des {\it isom\'etries globales}
     pour la m\'etrique
    sur $M$.

      Dans le cas d'une forme quadratique non d\'eg\'en\'er\'ee g\'en\'erique il n'existe pas ``autant"
       d'isom\'etries globales . 
      En revanche, dans une carte locale de $M$, donn\'ee par une
      section locale de la projection $G \to G/\Gamma$, les translations \`a droite d\'efinissent des transformations
      locales qui agissent par isom\'etrie.  Les m\'etriques holomorphes ainsi construites sont
      donc {\it localement homog\`enes} (pour tous les couples de points $(m,m')$ dans $M$, il existe une isom\'etrie 
      locale
      d'un voisinage ouvert de $m$ dans $M$ sur un voisinage ouvert de $m'$ qui envoit $m$ sur $m'$). Pour une forme quadratique non d\'eg\'en\'er\'ee g\'en\'erique
      la m\'etrique riemannienne holomorphe obtenue est un exemple de structure g\'eom\'etrique de type g\'en\'eral (voir la section suivante pour la d\'efinition) pour laquelle le corollaire \ref{corollaire} s'applique.     
      
      Remarquons aussi que si la m\'etrique riemannienne holomorphe qui
      correspond \`a la forme de Killing est de courbure sectionnelle constante, pour les autres
      exemples la courbure sectionnelle est, en g\'en\'eral, non constante. Elle est repr\'esent\'ee, en g\'en\'eral, par une fonction m\'eromorphe non constante d\'efinie
      sur la $2$-grasmannienne des plans non d\'eg\'en\'er\'es (qui constitue un ouvert de la $2$-grasmannienne de l'espace tangent holomorphe \`a $M$).
      
      Certains quotients par des groupes finis des exemples pr\'ec\'edents poss\`edent \'egale-ment des m\'etriques riemanniennes holomorphes. Pour s'y convaincre sur un exemple
      tr\`es simple, prenons un tore complexe de dimension $2$ de la forme $\CC^2$ quotient\'e par un r\'eseau $\Lambda$. Le quotient du tore par une transformation d'ordre
      $2$ du type $\displaystyle  (z_{1}, z_{2})  \to  (z_{1} + \frac{\lambda}{2}, -z_{2}),$  o\`u $\lambda \in \Lambda$ est un g\'en\'erateur du r\'eseau,  poss\`ede  la  m\'etrique riemannienne holomorphe induite par la m\'etrique plate $dz_{1}^2 + dz_{2}^2$, invariante par les translations et la transformation consid\'er\'ee.

     En dimension $3$, des exemples  in\'edits de m\'etriques riemanniennes holomorphes sur des vari\'et\'es complexes 
     dont le fibr\'e tangent n'est pas holomorphiquement trivial ont \'et\'e construits dans~\cite{Gh}. Ces exemples 
     s'obtiennent \`a
     partir des espaces homog\`enes du groupe $SL(2,\CC)$, par d\'eformation  de la structure complexe.
     
      L'id\'ee est de consid\'erer un r\'eseau co-compact $\Gamma$ de $SL(2,\CC)$ et de perturber l'action par 
      translations \`a droite
 de $\Gamma$ sur $SL(2,\CC)$ de mani\`ere que le quotient soit toujours une vari\'et\'e. L'auteur prouve qu'il existe
 des morphismes de groupes $u : \Gamma \to SL(2,\CC)$ tels que l'action \`a droite de $\Gamma$ sur $SL(2,\CC)$ donn\'ee
 par :
 $$(m,\gamma) \in SL(2,\CC) \times \Gamma \to u(\gamma^{-1})m\gamma \in SL(2,\CC)$$ 
 est libre et totalement discontinue. Le quotient est une vari\'et\'e complexe compacte $M(u,\Gamma)$ sur laquelle la
 m\'etrique de Killing de $SL(2,\CC)$ induit une m\'etrique holomorphe. Il existe des morphismes $u$ pour lesquels
 les vari\'et\'es $M(u,\Gamma)$ ne poss\`edent aucun champ de vecteurs holomorphe non nul. 
 Il est prouv\'e dans~\cite{Gh} que tout
 tenseur holomorphe sur une vari\'et\'e de type $M(u,\Gamma)$  se rel\`eve en un tenseur holomorphe sur $SL(2,\CC)$ invariant
 par l'action \`a droite de $SL(2,\CC)$ et par l'action \`a gauche de $u(\Gamma)$. 
 En particulier, tout tenseur holomorphe (et donc toute m\'etrique riemannienne holomorphe) sur $M(u,\Gamma)$ est localement homog\`ene.  Il s'agit d'un cas 
 particulier du th\'eor\`eme~\ref{homogene}.
  
\section{Structures g\'eom\'etriques rigides}

Dans cette section nous rappelons bri\`evement le concept de  structure g\'eom\'etrique (de type fini) introduit pour la premi\`ere fois par C. Ehresmann pour
donner une interpr\'etation g\'eom\'etrique globale aux  travaux de E. Cartan. Nous suivrons \'egalement la m\'ethode expos\'ee par M.Gromov dans \cite{Gro} pour comprendre les isom\'etries des
structures {\it rigides}.

 Consid\'erons $M$ une vari\'et\'e complexe de dimension $n$ et d\'esignons par $R^r(M)$ le fibr\'e des $r$-rep\`eres de $M$. Il s'agit d'un fibr\'e principal au-dessus
 de $M$ de groupe structural $D^r(\CC^n)$, le groupe des $r$-jets en $0$ de biholomorphismes locaux de $\CC^n$ qui pr\'eservent $0$. Par exemple,
 le fibr\'e des $1$-rep\`eres s'identifie avec l'ensemble des bases de tous les espaces tangents \`a $M$ et $D^1(\CC^n)$ n'est rien d'autre que $GL(n, \CC)$.
 
 \begin{definition}
 
 Une structure g\'eom\'etrique  holomorphe de type affine et d'ordre $r$  sur $M$
 est une application holomorphe $D^r(\CC^n)$-equivariante du fibr\'e $R^r(M)$ dans une vari\'et\'e affine $\cal A$  munie d'une action alg\'ebrique du groupe $D^r(\CC^n)$.
 
 \end{definition}
 
 Par un th\'eor\`eme de plongement classique d\^u  \`a Chevalley, il est toujours possible de consid\'erer que la vari\'et\'e affine $\cal A$ est plong\'ee de mani\`ere
 \'equivariante dans un espace vectoriel complexe de dimension finie muni d'une action alg\'ebrique de $D^r(\CC^n)$. Tout tenseur holomorphe sur $M$ est une structure
 g\'eom\'etrique holomorphe de type affine d'ordre $1$, tandis que les connexions lin\'eaires holomorphes sont des structures de type affine d'ordre $2$ :  le th\'eor\`eme~\ref{homogene}
 s'applique bien \`a ces structures g\'eom\'etriques. Un champ de droites holomorphe ou une structure conforme holomorphe sont des structures g\'eom\'etriques holomorphes qui ne sont pas de type affine (mais de type projectif) et le th\'eor\`eme~\ref{homogene} ne s'applique pas dans ce cadre. Pour un r\'esultat d'homog\'en\'eit\'e locale concernant les
 structures conformes holomorphes le lecteur pourra se r\'ef\'erer \`a~\cite{Be}.
 
    Le jet d'ordre $s$ d'une structure g\'eom\'etrique $\varphi$ d'ordre $r$ est une nouvelle structure g\'eom\'etrique
    d'ordre $s+r$. Si dans une carte locale de la vari\'et\'e $\varphi$ est d\'etermin\'ee par une application \`a valeurs dans
    la vari\'et\'e $\cal A$, son jet d'ordre $s$, que nous notons $\varphi^{(s)}$, s'exprime dans la m\^eme carte par le $s$-jet de
    l'application initiale. 
    
    En la pr\'esence d'une m\'etrique riemannienne holomorphe il convient de privil\'egier les syst\`emes de coordonn\'ees exponentielles et dans ce cas il n'est
    pas n\'ecessaire de consid\`erer le fibr\'e des $r$-rep\`eres. Rappelons \`a cet effet que chaque \'el\'ement du fibr\'e des $1$-rep\`eres $R^1(M)$ repr\`esentant
    par d\'efinition un $1$-jet de carte est r\'ealis\'e par une unique carte exponentielle pour la m\'etrique riemannienne holomorphe. Autrement dit, chaque
    \'el\'ement du fibr\'e $R^1(M)$ fournit une carte (exponentielle) de la vari\'et\'e $M$ et d'autant plus un \'el\'ement du fibr\'e $R^r(M)$.  Gr\^ace \`a cette
    section du fibr\'e des $1$-rep\`eres \`a valeurs dans le fibr\'e de $r$-rep\`eres, on pourra consid\'erer que toutes les structures g\'eom\'etriques
    sont d\'efinies sur $R^1(M)$.
    
   La pr\'esence de la m\'etrique riemannienne holomorphe $q$ nous permet de consid\'erer le fibr\'e des $1$-rep\`eres $q$-orthonorm\'es et, quitte
   \`a prendre un rev\^etement double non ramifi\'e de $M$ (ce qui ne change rien \`a notre probl\`eme d'homog\'en\'eit\'e locale), de consid\'erer le fibr\'e des
   $1$-rep\`eres $q$-orthonorm\'es directs $R(M)$ qui est un sous-fibr\'e principal de $R^1(M)$ de groupe structural $SO(3,\CC)$.
   
   Par exemple, on peut  voir le $s$-jet $q^{(s)}$ d'une m\'etrique riemanninne holomorphe $q$ comme \'etant une application \'equivariante du $SO(3, \CC)$-fibr\'e des
   rep\`eres orthonorm\'es directs \`a valeurs dans l'espace affine $J^s$ 
         constitu\'e par les $s$-jets en $0$ de m\'etriques riemanniennes holomorphes
   sur $\CC^3$ dont le $1$-jet est $dz_1^2+dz_2^2+dz_3^2$. 
  L'espace affine   $J^s$  
  est muni de l'action
  alg\'ebrique  affine de $SO(3,\CC)$ obtenue par changement de carte exponentielle. Comme cette action fixe l'\'el\'ement $dz_1^2+dz_2^2+dz_3^2$  nous pouvons consid\'erer
  que $J^s$ est un espace vectoriel d'origine $dz_1^2+dz_2^2+dz_3^2$ muni d'une action lin\'eaire et alg\'ebrique de $SO(3, \CC)$.
  
  Une {\it isom\'etrie locale}  d'une structure g\'eom\'etrique $\varphi$ est un biholomorphisme local entre deux ouverts de $M$ qui pr\'eserve $\varphi$. L'ensemble $Is^{loc}$
  des isom\'etrie locales   est un pseudo-groupe pour la composition. Quand le pseudo-groupe des isom\'etries locales agit transitivement sur $M$, la structure
  g\'eom\'etrique $\varphi$ est dite {\it localement homog\`ene}. Un champ de Killing (local) pour $\varphi$ est un champ de vecteurs holomorphe (local) sur $M$ dont
  le flot (local ) agit par isom\'etries locales pour $\varphi$. Deux points suffisamment proches de $M$ qui sont dans la m\^eme orbite de $Is^{loc}$ sont reli\'es par
  une isom\'etrie locale proche de l'identit\'e et se trouvent donc dans la m\^eme orbite d'un champ de Killing local~\cite{Ben}.

 Un  th\'eor\`eme d\'emontr\'e par I.Singer~\cite{Si} pour les m\'etriques riemanniennes et g\'en\'eralis\'e dans \cite{Amo}, \cite{Gro}, \cite{D1} pour des structures g\'eom\'etriques plus g\'en\'erales affirme qu'il existe un $s$ suffisament grand tel que  $q$ est localement homog\`ene si et seulement si l'image de l'application $q^{(s)}$  est exactement 
  une $SO(3,\CC)$-orbite de $J^s$.  Plus pr\'ecis\'ement, il existe un entier positif $s$ (qui ne d\'epend que de $M$ et de $q$) tel que deux points $m$ et $m'$ sont reli\'es par une isom\'etrie locale de $q$ si et seulement si $q$-admet un m\^eme $s$-jet aux points $m$ et $m'$. De plus (voir, par exemple, la jolie preuve de \cite{DG}, {\it toute application lin\'eaire qui envoie une base orthonorm\'ee directe de $(T_{m}M, q_{m})$ sur une base orthonorm\'ee directe de $(T_{m'}M, q_{m'})$  pour laquelle le $s$-jet de $q$ est le m\^eme
  d\'efinit en coordonn\'ees exponentielles un germe de biholomorphisme qui envoie $m$ sur $m'$ et qui est une isom\'etrie locale}.
  
  Comme il est montr\'e dans \cite{Gro} et \cite{DG} le ph\'enom\`ene pr\'ec\'edent est caract\'eristique aux structures g\'eom\'etriques que M. Gromov appelle
  {\it rigides} et  qui sont caract\'eris\'ees par la propri\'et\'e que le pseudo-groupe des isom\'etries locales est de dimension finie. Le fait qu'une  m\'etrique riemannienne holomorphe
  soit  une  structures rigide est  essentiellement  d\^u au fait qu'une isom\'etrie locale est compl\'etement d\'etermin\'ee par son jet d'ordre $1$ en un point.
  
  Par la suite, quelque soit la structure g\'eom\'etrique (rigide ou non) qui existe sur $M$ (penser \`a un tenseur), nous allons consid\'erer la structure plus riche $\varphi$ 
  qui est la juxtaposition de la structure initiale  et de la m\'etrique riemannienne holomorphe $q$. Comme $q$ est rigide, $\varphi$ sera (d'autant plus) une structure g\'eom\'etrique
  rigide et on aura encore  que {\it toute application lin\'eaire qui envoie une base orthonorm\'ee directe de $(T_{m}M, q_{m})$ sur une base orthonorm\'ee directe de $(T_{m'}M, q_{m'})$  pour laquelle le $s$-jet de $\varphi$ est le m\^eme,
  d\'efinit en $q$-coordonn\'ees exponentielles un germe de biholomorphisme qui envoie $m$ sur $m'$ et qui est une isom\'etrie locale pour $\varphi$}. Ce proc\'ed\'e
  ne prend pas en compte  toutes les isom\'etries locales de la structure initiale, mais seulement celles qui pr\'eservent en m\^eme temps la m\'etrique riemannienne
  holomorphe $q$. Notre m\'ethode vise donc \`a d\'emontrer directement que  la structure g\'eom\'etrique plus riche qui englobe \`a la fois la structure initiale et $q$ est encore  localement homog\`ene, ce qui implique en particulier l'homog\'en\'eit\'e locale de la structure initiale.
  
 \begin{definition}  Une structure g\'eom\'etrique d'ordre $r$ et de type affine $\varphi$ sur $M$ sera dite de type g\'en\'eral s'il existe des  entiers $s$ arbitrairement grands tels que
   l'image de l'application $s$-jet de
 $\varphi$ contiennent au moins une orbite de stabilisateur fini dans $D^{r+s}(\CC^n).$
 
 \end{definition}
  
  Pour une structure rigide ceci veut dire qu'il existe au moins un point $m \in M$ au voisinage duquel il n'existe pas de champ de Killing local s'annulant en $m$. En effet,
  la rigidit\'e implique que pour un $s$ suffisamment grand le stabilisateur du $s$-jet de $\varphi$ en un point $m$ se prolonge de mani\`ere unique en un \'el\'ement du pseudo-groupe des
  isom\'etries locales qui fixe $m$. Un sous-groupe \`a un param\`etre du stabilisateur du $s$-jet de $\varphi$ en $m$ fournit donc un sous-groupe \`a un param\`etre d'isom\'etries locales
  pour $\varphi$ qui fixent $M$. Par un th\'eor\`eme classique d\^u \`a S. Lie un tel sous-groupe \`a un param\`etre d'isom\'etries n'est rien d'autre que le flot local d'un champ de Killing s'annulant en $m$.
  
  Finissons cette section par la preuve du fait que le th\'eor\`eme \ref{homogene} implique le corollaire \ref{corollaire}.
  
  Consid\'erons une vari\'et\'e complexe compacte $M$ de dimension $3$ munie d'une m\'etrique riemannienne holomorphe $q$ et d'une structure g\'eom\'etrique de type affine et de type g\'en\'eral (cette structure g\'eom\'etrique pourrait tr\`es bien \^etre la m\'etrique $q$ elle-m\^eme : on a vu de tels exemples de m\'etriques \`a la section $2$).
  D'apr\`es le th\'eor\`eme \ref{homogene} la structure g\'eom\'etrique $\varphi$ form\'ee par la juxtaposition de $q$ et de la structure g\'eom\'etrique de type affine et de type
  g\'en\'eral est une structure g\'eom\'etrique rigide de type affine et de type g\'en\'eral qui est {\it localement homog\`ene}.  Ceci equivaut \`a~:  $\forall s \in \NN$ l'image
  de l'application $\varphi^{(s)}$, d\'efinie sur le fibr\'e des rep\`eres $q$-orthonorm\'ees directes $R(M)$, est constitu\'ee d'une unique orbite qui s'identifie au quotient de $SO(3, \CC)$ par le stabilisateur d'un point. Nous avons expliqu\'e un peu plus haut pourquoi il convient de voir $SO(3, \CC)$ comme \'etant plong\'e canoniquement (par le choix de $q$) dans
  $D^{r+s}(\CC^3)$. Ceci montre que le stabilisateur du $s$-jet de $\varphi$ dans $SO(3, \CC)$ est encore un sous-groupe fini $F$. Nous avons donc une application
  holomorphe $SO(3, \CC)$ \'equivariante du fibr\'e $R(M)$ dans le quotient $SO(3, \CC)/F$. Il est \'equivalent de dire qu'il existe un rev\^etement fini non ramifi\'e de $M$
  associ\'e au groupe $F$ au-dessus duquel le fibr\'e principal $R(M)$ se trivialise. Ce rev\^etement fini de $M$ poss\`ede un fibr\'e  tangent holomorphe  trivial et d'apr\`es
  le th\'eor\`eme de Wang \cite{Wa} il s'identifie \`a un quotient d'un groupe de Lie complexe connexe et simplement connexe de dimension $3$ par un r\'eseau co-compact.

\section{Homog\'en\'eit\'e locale sur un ouvert dense}

Soit $(M,q)$ une vari\'et\'e complexe compacte connexe de dimension $3$ munie d'une m\'etrique riemannienne holomorphe. 
Nous d\'emontrons dans cette section la version faible suivante du th\'eor\`eme \ref{homogene}  : {\it toute structure g\'eom\'etrique de type affine $\psi$ sur $M$ est localement homog\`ene
sur un ouvert dense (en dehors d'un sous-ensemble analytique compact, eventuellement vide)}. Il est donc \'egalement vrai que la structure  g\'eom\'etrique $\varphi =(\psi, q)$ obtenue  
par la juxtaposition de $\psi$ et $q$ est localement homog\`ene sur un ouvert dense. En particulier la m\'etrique riemannienne holomorphe $q$ est localement homog\`ene sur un ouvert dense. Ce r\'esultat
pr\'esente d\'ej\`a un progr\`es notable par rapport au th\'eor\`eme obtenu dans \cite{D2} o\`u on avait seulement pu d\'emontrer que les orbites du pseudo-groupe des isom\'etries
locales forment ou bien un ouvert dense ou bien sont les fibres d'une fibration sur une courbe de genre $g \geq  2$. Le travail fait dans cette section \'elimine la deuxi\`eme
alternative en pr\'ecisant la structure de la fibration $\pi$. Remarquons pour finir que dans~\cite{D2} les r\'esultats sont \'enonc\'es seulement pour la m\'etrique riemannienne
holomorphe, mais il suffit  de remplacer syst\`ematiquement dans la preuve le ''$s$-jet de $q$'' avec le  ''$s$-jet de la structure g\'eom\'etrique $\phi = (\psi, q)$'' pour obtenir les m\^emes r\'esultats  avec $\phi= (\psi, q)$ \`a la place de $q$.

\subsection{ $M$ est un fibr\'e principal en tores complexes}

Nous utilisons ici de mani\`ere essentielle le lemme 6.4 que nous avons prouv\'e dans \cite{D2}. Il s'agit d'un lemme qui exploite le manque d'homog\'en\'eit\'e locale
pour fabriquer certains tenseurs holomorphes :

 \begin{lemme}  \label{construction}
 
 Soit $M$ une vari\'et\'e complexe compacte de dimension $3$ munie d'une m\'etrique riemannienne
 holomorphe $q$ et \'egalement d'une structure g\'eom\'etrique de type affine $\psi$  (il se peut que $\psi$=  $q$).
 
 i). Si la structure g\'eom\'etrique $\varphi =(\psi,q)$, form\'ee par la juxtaposition de $\psi$ et de $q$, n'est pas localement homog\`ene (en particulier, si $q$ ou $\psi$ ne sont pas localement homog\`enes) sur $M$, alors il existe un entier
 positif $n$ et
 une section holomorphe non triviale du fibr\'e vectoriel $S^n(T^{*}M)$
  qui s'annule en au moins un point et 
 qui prend
 ses valeurs dans le c\^one des puissances $n$-i\`emes.
 
 ii). Si  la structure g\'eom\'etrique $\varphi =(\psi, q)$  n'est localement homog\`ene sur aucun ouvert dense de $M$, alors il existe  
   un entier positif $n$ et deux sections holomorphes  lin\'eairement 
 ind\'ependantes du
 fibr\'e $S^n(T^{*}M)$ qui s'annulent en au moins un point et qui prennent leurs valeurs dans le c\^one des
 puissances $n$-i\`emes. Ces sections sont en tout point colin\'eaires.
 
\end{lemme}

Supposons par l'absurde  qu'il existe une structure g\'eom\'etrique holomorphe de type affine $\psi$  sur $M$ telle que $\varphi=(\psi, q)$ n'est localement homog\`ene sur aucun ouvert dense de $M$.  Le point ii) du lemme~\ref{construction} ensemble avec le lemme suivant qui est une cons\'equence directe d'un th\'eor\`eme de M. Brunella~\cite{Bru} nous
montre que $M$ fibre sur une courbe de genre $g \geq 2$. Le lemme est le suivant :

\begin{lemme} \label{fibration}
 Soit $M$ une vari\'et\'e complexe compacte de dimension $3$ avec  fibr\'e canonique  trivial. Supposons qu'il existe  un entier positif $n$ et deux sections holomorphes 
 lin\'eairement ind\'ependantes $\alpha_1$ et   
 $\alpha_2$  du fibr\'e  $S^n(T^{*}M)$ qui prennent leurs valeurs dans le c\^one des puissances $n$-i\`emes et qui
 sont en tout point colin\'eaires.
 
 Il existe alors une fibration holomorphe $\pi$ de $M$ sur une surface de Riemann $C$ de genre $g \geq 2$ et deux
 sections holomorphes $\eta_1$ et $\eta_2$ de $(T^{*}C)^{\otimes n}$ telles que $\alpha_i=\pi^{*}(\eta_i)$.
\end{lemme}

 \begin{demonstration}
 
Les arguments de~\cite{Bru}, notamment le lemme $1$ de la page $4$, la proposition $3$ de la section $3$   et  la preuve du cas  $1$ de la section $4$ (page $15$)
s'appliquent  sans aucune modification d\`es qu'on d\'emontre que les formes $\alpha_1$ et   
$\alpha_2$ sont int\'egrables. L'int\'egrabilit\'e a \'et\'e montr\'e dans  la preuve du lemme 6.6 de~\cite{D2}.
\end{demonstration}

Dans~\cite{D2} nous avions montr\'e que les fibres de $\pi$ sont form\'ee par des points qui sont dans la m\^eme orbite
du pseudo-groupe des isom\'etries locales de $\varphi =(\psi, q)$.  La strat\'egie \`a suivre sera maintenant la suivante : nous \'etudions cette fibration et  d\'emontrons  que
$M$ admet une structure de fibr\'e principal en tores complexes (avec eventuellement des fibres singuli\`eres) au-dessus de la courbe $C$. Finalement un lemme de prolongement de
formes diff\'erentielles multiples inspir\'e par~\cite{Ue} permet d'obtenir une contradiction et de conclure que $\varphi$ doit \^etre localement homog\`ene sur un ouvert dense.

 Regardons de plus pr\`es la structure de la fibration pr\'ec\'edente et d\'emontrons d'abord que les fibres
g\'en\'eriques sont des tores complexes.  On prouve ensuite   que ces tores sont n\'ecessairement isomorphes, ce qui donne :

\begin{proposition}

 $M$ est un fibr\'e principal, de fibre type  un  tore complexe de dimension  $2$, au-dessus
d'une surface de Riemann de genre $g \geq 2$ (avec eventuellement un nombre fini de fibres singuli\`eres).

\end{proposition}

\begin{demonstration}

Les fibres de la fibration $\pi$  sont les feuilles du feuilletage donn\'e par le noyau de la section $\alpha_{1}$.
On peut d\'ecrire la distribution $Ker \alpha_{1}$ comme \'etant  \'egalement l'orthogonal $Y_{1}^{\bot}$ (au sens de $q$)
de la section $Y_{1}$ de $S^n(TM)$ duale \`a $\alpha_{1}$.  

Rappelons que la norme de $Y_{1}$ (au sens de $q^{\otimes n}$) est une fonction holomorphe et, par
cons\'equent, constante sur la vari\'et\'e compacte $M$. Cette constante est n\'ecessairement nulle
car $Y_{1}$ admet au moins un point d'annulation et la section $Y_{1}$ est donc en tout point isotrope.
La restriction de $Y_{1}$ \`a une fibre r\'eguli\`ere de $\pi$ sur laquelle $Y_{1}$ ne s'annule pas (rappelons que
$Y_{1}$ est projetable sur la base) munit cette fibre d'un champ de vecteurs holomorphe tangent non singulier et $q$-isotrope. Ces fibres sont donc n\'ecessairement des surfaces minimales : le flot de ce champ de vecteurs non singulier 
permet de ``bouger''  toute eventuelle copie de $P^1(\CC)$ incluse dans la fibre, ce qui implique que l'auto-intersection de $P^1(\CC)$  ne peut pas \^ etre n\'egative. 

De mani\`ere duale la restriction de $\alpha_1$ permet de construire le long d'une fibre r\'eguli\`ere
de la fibration une $1$-forme diff\'erentielle holomorphe $\alpha$ dont le noyau co\"{\i}ncide avec
l'espace tangent \`a la fibre.
Les fibres r\'eguli\`eres  poss\`edent  une forme
volume holomorphe (leur fibr\'e canonique est trivial) : cette forme volume est donn\'ee, par exemple, par la formule
$vol_{fibre}  (x,y)=\frac{vol(x,y,z)}{\alpha(z)},$ o\`u $vol$ d\'esigne la forme volume de $M$, $x,y$ des vecteurs
tangents \`a la fibre et $z$ n'importe quel vecteur transverse \`a la fibration  (il est ais\'e de voir que l'expression du volume de la fibre ne d\'epend pas
du suppl\'ementaire $z$ choisi).
La classification des surfaces complexes (voir, par exemple,~\cite{BCV}) montre que la fibre g\'en\'erique est ou bien un tore complexe
ou bien une surface de Kodaira primaire (une fibration non singuli\`ere en courbes elliptiques sur une courbe elliptique). Ce dernier cas est
\'elimin\'e par  le  th\'eor\`eme 4.1 de~\cite{Ue} (il s'agit d'un cas facile de ce th\'eor\`eme).

La prochaine \'etape montre que ces tores sont holomorphiquement isomorphes.  Soit $m$
un point de la base $C$ qui est une valeure r\'eguli\`ere de la fibration $\pi $ et soit $T_{m} = \pi^{-1} (m)$ la  fibre r\'eguli\`ere correspondante. La section isotrope  $Y_{1}$ permet de
construire  en restriction  au tore complexe $T_{m}$ un champ de vecteurs constant $Y$ et on peut
d\'efinir le long de $T_{m}$ un deuxi\`eme champ de vecteurs $X$ avec la propri\'et\'e que $q(X)=1$ et $X \in Y^{\bot}$ (on utilise ici que la fibre g\'en\'erique est un tore et donc son fibr\'e tangent est trivial : il suffit donc de prendre un deuxi\`eme champ de vecteurs tangent \`a la fibre et non colin\'eaire \`a $Y$ et de le diviser par sa $q$-norme qui est constante gr\^ace
\`a la compacit\'e de la fibre).  Compl\'etons 
ces deux champs de vecteurs par un troisi\`eme champ $Z$ en imposant la condition $q(Z)=0$, $q(Z,X)=0$ et $q(Y,Z)=1$. Cette condition exprime
le fait que $Z$ engendrent la deuxi\`eme droite isotrope du champ de plans non deg\'en\'er\'es $X^{\bot}$ et que $Z$ est choisi sur cette droite de 
mani\`ere unique en imposant  que la base $(X,Y,Z)$ soit de volume $1$. Ces trois champs constituent alors une trivialisation  du fibr\'e des rep\`eres
orthonorm\'es au-dessus de $T_{m}$. Ceci implique que le $1$-jet de biholomorphisme qui envoie un point $s \in T_{m}$ en un point  $s' \in T_{m}$ 
et dont la diff\'erentielle envoie $(X(s), Y(s), Z(s))$ sur $(X(s'), Y(s'), Z(s'))$ s'int\`egre en une isom\'etrie locale (le $s$-jet de la structure g\'eom\'etrique $(q, \varphi)$ est le m\^eme en les bases
$(X(s), Y(s), Z(s))$ et $(X(s'), Y(s'), Z(s'))$ \`a cause de la compacit\'e de la fibre $T_{m}$ et du fait que le fibr\'e des rep\`eres est trivial au-dessus de $T_{m}$). Cette isom\'etrie locale envoie
n\'ecessairement le temps $t$ de la g\'eod\'esique issue du point $s$ dans la direction $Z(s)$ sur le temps $t$ de la g\'eod\'esique issue du point $s'$ dans la direction $Z(s')$:
autrement dit, pour $t$ suffisament proche de $0$ l'image du point $exp_{s}(tZ(s))$ est $exp_{s'}(tZ(s'))$. Le temps $t$ du  flot g\'eod\'esique du champ $Z$ (d\'efini le long
de la fibre en $m$) envoie donc la fibre en $m$ en un ensemble de points reli\'es par des isom\'etries locales. Cet ensemble est contenu  donc dans la m\^eme fibre de la fibration
 $\pi$~:
rappelons que la conclusion du th\'eor\`eme principal de~\cite{D2} est pr\'ecis\'ement que les fibres de la fibration $\pi$ sont les orbites du pseudo-groupe des
isom\'etries locales de $\varphi$.
Le flot g\'eod\'esique de $Z$ r\'ealise alors un isomorphisme entre  la fibre en $m$ et les fibres proches. Les fibres r\'eguli\`eres de la fibration $\pi$ sont donc isomorphes
entre elles.
 
Comme un biholomorphisme proche de l'identit\'e d'un tore complexe est une translation, le flot 
g\'eod\'esique du champ transverse $Z$ identifie  les fibres voisines de la fibration par un biholomorphisme qui est une translation. 
Il vient qu'en dehors des fibres singuli\`eres la fibration $\pi$ admet une structure de fibr\'e principal.   
\end{demonstration}

\subsection {Prolongement de certaines sections de $K^{\otimes n}$}

Nous venons de voir que $M$ admet une structure fibr\'ee principale de fibre un  tore complexe $T$ de dimension $2$  au-dessus d'une courbe de genre $g \geq 2$, avec eventuellement des
fibres singuli\`eres. Si $u,v$ sont des coordonn\'ees
complexes sur la fibre type $T$ de la fibration en question, alors la $2$-forme diff\'erentielle holomorphe $\theta = du \wedge dv$ est invariante par les translations de la fibre
et est donc bien d\'efinie sur $M$ en dehors des fibres singuli\`eres.  Aussi  les sections $\omega_{1}=\alpha_{1}  \wedge \theta ^{\otimes n}$ et 
$\omega_{2}=\alpha_{2}  \wedge \theta^{\otimes n}$ 
de $K^{\otimes n}$, o\`u $K$ est le fibr\'e canonique de $M$,  sont bien d\'efinies en dehors des fibres singuli\`eres et non identiquement nulles. Nous allons arriver \`a une contradiction en  prolongeant les sections $\omega_{1}$ et
$\omega_{2}$ \`a l'aide du lemme suivant du \`a Ueno~\cite{Ue} :

\begin{lemme}   \label{prolongement}

Soit $M$ une vari\'et\'e complexe, $S$  un sous-ensemble analytique
compact de $M$ et $\omega$ une section holomorphe au-dessus de $M \setminus S$ de  $K^{\otimes n}$, o\`u $K$ est le fibr\'e canonique de $M$. Si $\int_{M \setminus S}  {| \omega  | ^{\frac{2}{n} }}  < \infty , $ alors $\omega$ est une section m\'eromorphe de $K^{\otimes n}$  admettant un p\^ole
d'ordre au plus $n-1$ en $S$.

\end{lemme}

Dans le lemme pr\'ec\'edent $| \omega | ^{\frac{2}{n} } $ d\'esigne la $2m$-forme diff\'erentielle r\'eelle qui s'exprime localement  par 
$|f(z) | ^{\frac {2}{n}} (dx_{1} \wedge dy_{1} \wedge \ldots \wedge dx_{m} \wedge dy_{m})$, o\`u $z_{k}=x_{k} + iy_{k}$, $m$ est la dimension complexe de $M$ et  $\omega = f(z)(dz_{1} \wedge  \ldots  \wedge dz_{m})^n$
est l'expression locale de $\omega$.

Nous allons donc appliquer  le lemme~\ref{prolongement}  \`a notre situation pour prolonger les sections $\omega_{1}$ et $\omega_{2}$ au-dessus des fibres singuli\`eres. Les arguments
suivants s'inspirent de~\cite{Ue}. 

Soit $c \in C$ une valeure critique  de notre  fibration principale,  $\pi ^{-1}(c)$ la fibre singuli\`ere correspondante
et $t$ une coordonn\'ee locale sur $C$ au voisinage de $c$ qui s'annule en $c$. Consid\'erons un disque $D= \{t, |t | < \epsilon \}$ suffisamment petit pour que $\pi ^{-1}(0)$ soit l'unique fibre singuli\`ere
contenue dans $\pi ^{-1} (D)$. 
Rappelons que $\alpha_1=\pi^{*}(\eta_1)$, o\`u $\eta_1$ est une section holomorphe d'une puissance du  fibr\'e canonique de $C$. Choisissons la coordonn\'ee locale
$t$ telle que $\eta_{1}(t) = t^l (dt) ^{\otimes n} $, pour un certain entier positif $l$ (qui  vaut $0$ dans le cas o\`u $\eta_{1}$ est non singuli\`ere en $c$).
  Dans ce cas l'expression locale de $\omega_{1}$, en restriction \`a l'ouvert $\pi ^{-1}(D \setminus \{0\})$ est $t^l (dt \wedge du \wedge dv)^n$. Pour montrer que
  $\omega_{1}$ se prolonge en une section holomorphe sur $\pi ^{-1}(D)$, il suffit de prouver que $(dt \wedge du \wedge dv)^n$ se prolonge. Pour cela on applique
  le lemme~\ref{prolongement} \`a la section $\omega = t^{-(n-1)} (dt \wedge du \wedge dv)^n$. Notons d'abord $A= \int_{T} \theta \wedge {\bar \theta} $, le volume 
  de la fibre r\'eguli\`ere et appliquons le th\'eor\`eme de Fubini pour estimer l'int\'egrale suivante~:
  
 $$ \int_{\pi ^{-1}(D \setminus \{0\})} { | \omega | ^{\frac{2}{n}}}=A \cdot | \int_{D \setminus \{0\}} {t^{\frac{-2(n-1)}{n}} dt \wedge {\bar dt}}| = 4 \pi A \cdot \int_{0}^{\epsilon} {r^{-(1- \frac{1}{n})}  dr  } <  \infty.$$
  Le lemme~\ref{prolongement} s'applique et montre que la section $\omega = t^{-(n-1)} (dt \wedge du \wedge dv)^n$ est m\'eromorphe sur $\pi^{-1}(D)$ et admet un p\^ole
  d'ordre au plus $n-1$ sur la fibre singuli\`ere. Il s'ensuit donc que la section $(dt \wedge du \wedge dv)^n$ est holomorphe sur $\pi^{-1}(D)$ et que $\omega_{1}$ se
  prolonge au-dessus de la fibre singuli\`ere $\pi ^{-1}(c)$. Le m\^eme raisonnement s'applique aussi \`a $\omega_{2}$ et montre que  $\omega_{1}$ et $\omega_{2}$
  d\'efinissent deux sections holomorphes globales du fibr\'e $K^{\otimes n}$. Comme $\alpha_{1}= f  \cdot \alpha_{2}$ pour une certaine fonction m\'eromorphe non constante $f$, ceci
  reste vrai pour les sections $\omega_{1}$ et $\omega_{2}$. La contradiction vient du fait que le fibr\'e (trivial) $K^{\otimes n}$  ne peut pas admettre deux sections holomorphes
  globales dont le quotient est une fonction m\'eromorphe non constante sur $M$.

\section{\'Etendre l'ouvert dense \`a $M$}

Avec les m\^emes notations que dans la section pr\'ec\'edente,
il s'agit de prouver ici  que l'ouvert dense d'homog\'en\'eit\'e locale est toujours \'egal \`a $M$ tout entier, ce qui finit la preuve du th\'eor\`eme \ref{homogene}.

  La propri\'et\'e suivante de prolongement  de champs de Killing, d\'ecouverte pour la premi\`ere fois par K. Nomizu dans le cadre des m\'etriques riemanniennes analytiques~\cite{No} et g\'en\'eralis\'ee dans~\cite{Amo} et~\cite{Gro} pour les structures g\'eom\'etriques rigides analytiques est essentielle pour la suite : 
   chaque point $m \in M$ admet un voisinage ouvert $U_m$ avec la propri\'et\'e que tout champ
             de Killing holomorphe, d\'efini sur un ouvert connexe $U \subset U_m$, se prolonge en un champ de Killing
             sur $U_m$.  Autrement dit, {\it la fibre du faisceau des germes de champs de Killing d'une structure holomorphe rigide est une  alg\`ebre de Lie $\cal G$ qui ne d\'epend pas
             du point. }Dans notre cas cette alg\`ebre est de dimension complexe au moins $3$ car elle est transitive sur $M \setminus S$. Par ailleurs, comme l'action de $\cal G$ est
             libre sur le fibr\'e des rep\`eres qui est de dimension $6$, la dimension de $\cal G$ est inf\'erieure ou \'egale \`a $6$. Il reste que la dimension de $\cal G$
             est \'egale \`a $3$, $4$ ou $5$.

Supposons par l'absurde qu'il existe sur $M$ une certaine structure g\'eom\'etrique holomorphe $\psi$ qui  n'est pas localement homog\`ene partout, mais seulement sur un ouvert dense. 
Autrement dit, si $s$ est un entier positif suffisamment grand l'image de l'application $\psi ^{(s)}$  qui repr\'esente le $s$-jet de $\psi$ (et dont
les fibres  se projetent sur $M$ en les orbites du pseudo-groupe des isom\'etries locales de $\psi$ ) est incluse dans l'adh\'erence d'une orbite
$O$ sous l'action du groupe $SO(3, \CC)$. Rappelons qu'un r\'esultat classique (voir, par exemple, ~\cite{Hum}, ~\cite{Mum}) affirme que pour les repr\'esentations alg\'ebriques  chaque orbite est ouverte dans son 
adh\'erence : 
{\it l'image de $\psi ^{(s)}$ est donc constitu\'e\'e de l'orbite $O$ \`a laquelle s'ajoute eventuellement  des orbites de dimension strictement inf\'erieure contenues dans l'adh\'erence
de $O$.} 

D\'esignons par $S$ l'ensemble analytique compact form\'e par les points de $M$ o\`u le $s$-jet de $\psi$ appartient \`a $\bar O \setminus O$. Autrement dit, $M \setminus  S$ est l'ouvert maximal  de $M$ sur lequel le pseudo-groupe des isom\'etries locales de $\psi$ agit transitivement.

La preuve du  point i) du lemme~\ref{construction} (voir~\cite{D2} ) construit dans ce cas  une section holomorphe non triviale $Y$ d'un  fibr\'e vectoriel $S^n(T^{*}M)$ dont
le lieu d'annulation est  $S$ et
 qui prend ses valeurs dans le c\^one des puissances $n$-i\`emes. Il est utile pour la suite de rappeler que la pr\'esence de la m\'etrique riemannienne holomorphe $q$ d\'efinit un isomorphisme entre $T^{*}M$ et $TM$ et cette dualit\'e permet de voir $Y$ aussi comme une section de $S^n(TM).$

 Nous d\'emontrons la proposition suivante qui nous permettra de remplacer pour la suite la structure g\'eom\'etrique (quelconque) $\psi$, avec la structure g\'eom\'etrique $\varphi=(q,Y)$, form\'ee par la juxtaposition de $q$ et du tenseur  $Y$.

 \begin{proposition}
 
  L'ouvert $M \setminus S$ est  l'ouvert maximal de $M$ sur lequel le  pseudo-groupe des isom\'etries locales de la structure g\'eom\'etrique holomorphe $\varphi=(q,Y)$ (form\'ee par la juxtaposition de $q$ et de $Y$) agit transitivement.
 \end{proposition}
 
 \begin{demonstration}
 
Nous savons par la section pr\'ec\'edente que la structure g\'eom\'etrique $\varphi=(q, Y)$ est localement
homog\`ene sur un ouvert dense de $M$ (car d\'emontr\'e pour toutes les structures g\'eom\'etriques de type affine). 

L'ouvert dense sur lequel $\varphi$ est localement homog\`ene ne peut contenir
 des points o\`u $Y$ s'annule car dans ce cas $Y$ serait nul sur tout l'ouvert donc partout. Inversement, s'il existe au moins un point dans $M \setminus S$ o\`u le $s$-jet
 de $\varphi$ se trouve dans une orbite qui n'est pas dense dans l'image de $\varphi^{(s)}$, la preuve du lemme 6.4 de~\cite{D2} construit une section holomorphe non triviale d'un fibr\'e $S^n(T^{*}M)$
 qui s'annule au moins au point consid\'er\'e. Comme $Y$ ne s'annule pas en ce point, cette nouvelle section et $Y$ seront lin\'eairement ind\'ependantes et on est dans les hypoth\`eses d'application du lemme~\ref{fibration}. Le raisonnement
 de la section $4$ s'applique et fournit  l'existence d'une fibration de $M$ sur une courbe de genre $g \geq 2$ ce qui, comme on l'a vu, conduit \`a une  contradiction.  
 \end{demonstration}
 
Chaque orbite de dimension strictement inf\'erieure \`a $3$ qui se trouve dans l'image du
$s$-jet de $\varphi$  (par exemple, toute orbite contenue dans $\bar O \setminus O$) admet un stabilisateur de dimension sup\'erieure ou \'egale \`a $1$  dans $SO(3, \CC)$. Un tel stabilisateur fournit des isom\'etries fixant un point (autrement dit, des champs locaux de Killing pour $\varphi$ qui s'annulent en au moins un point et  qui se lin\'earisent donc en coordonn\'ees exponentielles).  Pr\'ecisons ceci : si $\varphi^{(s)}$
est l'application ''$s$-jet  de $\varphi$'' d\'efinie sur $R(M)$, alors l'image par $\varphi^{(s)}$ de la restriction de $R(M)$ au-dessus de $M \setminus S$ est exactement une orbite
$O$ de l'espace affine des $s$-jets de $\phi$, tandis que l'image de la restriction de $R(M)$ \`a $S$ est form\'ee par des orbites de dimension strictement inf\'erieures contenues
dans l'adh\'erence de $O$. Si $\varphi^{(s)}$ envoie la fibre de $R(M)$ au-dessus d'un point $u  \in S$ sur une  orbite de $\bar O \setminus O$ qu'on identifie  \`a $SO(3, \CC)$ quotient\'e
par un  stabilisateur $B$, alors la dimension complexe de $B$ est sup\'erieure ou \'egale \`a $1$ et pour $s$ suffisamment grand les \'el\'ements de $B$ s'int\`egrent
en des isom\'etries locales qui fixent le point $u$. Tout sous-groupe \`a un param\`etre de $B$ fournit donc un champ de Killing au voisinage de $u$ qui s'annule en $u$.
Un tel champ de Killing est donc lin\'earisable au voisinage de $u$ (car isom\'etrie de $q$). Pour la partie lin\'eaire d'un tel champ de Killing on a les deux possibilit\'es suivantes
(rappelons au passage   que  le groupe $SO(3, \CC)$ est isomorphe \`a $PSL(2, \CC)$ : pour s'en assurer il suffit de faire  agir $PSL(2, \CC)$ (par changement de variables) sur l'espace vectoriel (de dimension $3$) des formes quadratiques  \`a deux variables (de la forme $ax^2+bxy +cy^2$) en pr\'eservant le d\'eterminant) : 

\begin{enumerate}

\item  la  partie lin\'eaire est {\it semi-simple}, conjugu\'ee dans   l'alg\`ebre de Lie de $PSL(2,\CC)$  \`a un multiple de l'\'el\'ement 
$\left( 
   \begin{array}{cc}
                                                                 1    &  0  \\
                                                                 0     & -1 \\
                                                                 \end{array}
                                                                 \right)$  
     
    \item  la partie lin\'eaire du champ de Killing soit {\it unipotente}, autrement dit    conjugu\'ee dans $PSL(2, \CC)$ \`a l'\'el\'ement 
       $\left( 
   \begin{array}{cc}
                                                                 0   &   1 \\
                                                                 0     &  0 \\
                                                                 \end{array}
                                                                 \right)$.
                                                                 
  \end{enumerate}

Il est \'equivalent de dire que la partie lin\'eaire d'un champ de Killing est semi-simple  ou de pr\'eciser que la diff\'erentielle du flot du champ de Killing au point fixe stabilise
un champ de vecteurs de $q$-norme \'egale \`a $1$ dans l'espace tangent (forme quadratique  $xy$ dans l'isomorphisme avec $PSL(2, \CC)$). Aussi, il est \'equivalent de dire que la partie lin\'eaire d'un champ de Killing est unipotente   ou de pr\'eciser que la diff\'erentielle du flot de champ de Killing au point fixe stabilise un champ de vecteurs $q$-isotrope  dans l'espace tangent (forme quadratique  $x^2$). Dans ces situations on utilise
classiquement la terminologie : {\it isotropie semi-simple} et {\it isotropie unipotente}.

Le lemme suivant pr\'ecise  le type d'orbites qui peuvent se trouver dans l'image du $s$-jet de $\varphi$ et sera tr\`es utile  pour la suite :

\begin{lemme}  \label{dimension}

L'image de $\varphi^{(s)}$ contient une orbite dense $O$ de dimension $3$. Les orbites de $\bar O \setminus O$ contenues dans l'image de $\varphi^{(s)}$ sont de 
dimension $2$.
\end{lemme}

\begin{demonstration}

Commen\c cons par rappeler que  le seul  espace homog\`ene de dimension $1$ complexe de $PSL(2, \CC)$ est  la droite projective $P^1(\CC)$. En effet, le stabilisateur
d'un point doit \^etre un sous-groupe de Lie complexe de dimension $2$ de $PSL(2, \CC)$ et il est, par cons\'equent, conjugu\'e dans $PSL(2, \CC)$  au sous-groupe form\'e par les
matrices triangulaires sup\'erieures. Comme $P^1(\CC)$ est un espace compact, ce type d'orbite ne peut pas appara\^{\i}tre dans une repr\'esentation de $PSL(2, \CC)$
sur une vari\'et\'e affine. Dans l'image de $\varphi^{(s)}$ il n'existe donc aucune orbite de dimension $1$.

Montrons maintenant qu'il n'existe dans $\varphi^{(s)}$ aucune orbite de dimension $0$. Supposons par l'absurde le contraire. Il existe alors un point $u$ dans $S$ o\`u le stabilisateur
du $s$-jet de $q$ est de dimension (maximale) $3$. Consid\'erons l'alg\`ebre des germes des champs de Killing $\cal G$ au voisinage de ce point $u$. Le groupe d'isotropie en $u$ est de dimension $3$ isomorphe donc au groupe lin\'eaire $PSL(2, \CC)$. Le morphisme d'\'evaluation en $u$ qui \`a un champ de Killing associe sa valeur au point $u$
est un morphisme d'espaces vectoriels qui admet un noyau de dimension $3$. Comme la dimension de $\cal G$ est inf\'erieure ou \'egale \`a $5$, la dimension de l'image
de ce morphisme est un sous-espace vectoriel strict de $T_{s}M$ qui doit \^etre invariant par l'action lin\'eaire du groupe d'isotropie $PSL(2, \CC)$. Il vient que ce sous-espace
vectoriel est  n\'ecessairement trivial, ce qui implique que $\cal G$ est de dimension $3$ isomorphe \`a $PSL(2, \CC)$. La contradiction vient du fait qu'au voisinage
du point fixe $u$ les orbites de l'action lin\'earis\'ee (en coordonn\'ees exponentielles) de $PSL(2, \CC)$ sont  de dimension au plus $2$ (dans le mod\`ele lin\'eaire donn\'e par les formes quadratiques  \`a $2$ variables cette action pr\'eserve le d\'eterminant). L'action de $\cal G$  n'est donc pas  transitive en dehors de $S$ : absurde.

Il reste que les orbites de $\bar O \setminus O$ contenues dans l'image de $\varphi^{(s)}$ sont de dimension $2$. Ceci  implique que l'orbite dense $O$ est de dimension strictement sup\'erieure : elle est donc de dimension  $3$.
\end{demonstration}

La proposition suivante est une cons\'equence directe du lemme pr\'ec\'edent.

\begin{proposition}   \label{unimodulaire}

i) L'alg\`ebre des germes de champs de Killing  pour $\varphi$  est une alg\`ebre de Lie  non-unimodulaire de dimension $3$.

ii) L'ouvert  $M \setminus S$ admet une $(G,G)$-structure, o\`u $G$ est l'unique groupe de Lie connexe et simplement connexe associ\'e \`a $\cal G$. En restriction \`a $M \setminus S$,
$q$ provient d'une m\'etrique riemannienne holomorphe sur $G$ invariante par translations \`a gauche.

\end{proposition}

\begin{demonstration}

Au-dessus de $M \setminus S$ l'image du fibr\'e des rep\`eres dans l'espace des $s$-jets de $\varphi$ est l'orbite $O$  de dimension
  complexe $3$ qui s'identifie au quotient de $SO(3, \CC)$ par un sous-groupe fini. On b\'en\'eficie
  donc (eventuellement sur un rev\^etement fini non ramifi\'e) d'une trivialisation du fibr\'e
  des rep\`eres et, de plus, un biholomorphisme local qui relie deux points de $M \setminus S$ est une isom\'etrie locale
  si et seulement si son action pr\'eserve cette trivialisation .  
  
  L'action de $\cal G$ est donc libre et transitive sur $M \setminus S$ (isotropie triviale), ce qui implique que $\cal G$ est de dimension complexe $3$. 
  Avec le langage des $(G,X)-$ structures, nous sommes
  en pr\'esence d'une $(G,G)$-structure sur l'ouvert invariant $M \setminus S$, o\`u $G$ est l'unique  groupe de Lie complexe connexe simplement connexe de dimension $3$ associ\'e
  \`a l'alg\`ebre de Lie $\cal G$. Sur l'ouvert $M \setminus S$, la m\'etrique riemannienne holomorphe $q$ provient  d'une  m\'etrique riemannienne holomorphe sur $G$  invariante par les translations \`a gauche. On a vu dans la section $2$ qu'une telle m\'etrique  se construit en choisissant  une forme quadratique non d\'eg\'en\'er\'ee sur l'alg\`ebre de Lie $\cal G$
  et en la transportant par les translations \`a gauche. Comme le pseudo-groupe d'isom\'etries locales agit dans notre situation sans point fixe sur $M \setminus S$ (isotropie triviale), la forme quadratique $q$ n'est pr\'eserv\'ee par aucune transformation adjointe (voir~\cite{Mi}).

    Prouvons que  le groupe $G$ ne peut pas \^ etre    unimoduaire (il est, par cons\'equent, r\'esoluble)~\cite{Mi}. Si par l'absurde $G$ est  unimodulaire alors la transformation adjointe pr\'eserve le volume de la forme quadratique consid\'er\'ee sur l'alg\`ebre de Lie $\cal G$ ( qui fournit la restriction de $q$ \`a $M \setminus S$).  Soit  $(K_{1}, K_{2}, K_{3})$ une base fix\'ee de $\cal G$ et consid\'erons les
     trois champs de  Killing locaux lin\'eairement ind\'ependants d\'efinis au voisinage d'un point $m \in M \setminus S$ qui correspondent  \`a  la base $(K_{1}, K_{2}, K_{3})$ de $\cal G$. Si on agit  par une isom\'etrie locale qui envoie le point $m$ sur  un autre point $m'$ dans $M \setminus S$ les trois champs de Killing initiaux sont conjugu\'es \`a trois autres champs de Killing
     qui correspondent \`a l'image de la base $(K_{1},K_{2},K_{3})$ par la transformation adjointe (qui est suppos\'ee pr\'eserver le volume). L'action de $\cal G$ pr\'eserve
     donc le volume $vol(K_{1}, K_{2}, K_{3})$ des champs $K_{1}, K_{2}$ et $K_{3}$. Comme cette action est transitive ce volume est constant.
     
   Autrement dit, quelque soient
     trois champs de Killing locaux d\'efinis au voisinage d'un point de $M \setminus S$, leur volume par rapport \`a $q$ est une fonction constante. Pour obtenir la contradiction recherch\'ee consid\'erons un point $u$ dans $S$ et consid\'erons  un voisinage  $U$ de $u$ dans $S$ qui v\'erifie  la propri\'et\'e
     de prolongement de champs de Killing. Consid\'erons au voisinage d'un point de $U  \setminus S$ trois champs de Killing lin\'eairement ind\'ependants, correspondant  \`a une base $(K_{1}, K_{2}, K_{3})$ de $\cal G$. Ces champs de Killing se prolongent \`a $U$ en des champs de Killing dont les valeurs en $u$ sont n\'ec\'essairement des
     vecteurs de $T_{u}M$ lin\'eairement d\'ependants  (tangents \`a $S$ car $S$ est invariant sous l'action de  $\cal G$).  Le volume associ\'e \`a ces trois champs de
     Killing est donc une fonction constante  non nulle  sur $U  \setminus S$, qui se prolonge (holomorphiquement et donc continument) en une fonction qui s'annule en les points de $S$ : impossible.
\end{demonstration}

  \begin{proposition}
  
  Les  composantes connexes  de $S$ sont des  surfaces complexes lisses (sous-vari\'et\'es de codimension $1$ de $M$).
\end{proposition}

\begin{demonstration}
  
 Remarquons que   ${\bar O} \setminus O$ \'etant de dimension complexe $2$, les  orbites de dimension $2$ sont des ouverts de ${\bar O} \setminus O$. Par connexit\'e,  le fibr\'e des rep\`eres au-dessus de chaque composante connexe de $S$ est envoy\'e par l'application $s$-jet de $\varphi$ sur une m\^eme  orbite (de dimension 2) contenue dans  ${\bar O} \setminus O$ (et le stabilisateur $B$ d'une telle orbite dans $SO(3, \CC)$ est une extension par un groupe fini d'un sous-groupe \`a un param\`etre). Par cons\'equent le pseudo-groupe des isom\'etries locales de $\varphi = (q,Y)$  agit transitivement sur chaque composante connexe 
 de $S$. Il vient donc que chaque composante connexe de $S$ est {\it lisse}.  
 
 Au-dessus de chaque composante connexe de $S$ on b\'en\'eficie d'une application holomorphe $SO(3, \CC)$-\'equivariante du fibr\'e des rep\`eres $R(M)$ restreint \`a $S$
 dans l'espace homog\`ene $SO(3, \CC) /B$. Ceci s'interpr\`ete comme la donn\'ee d'un champ de vecteurs holomorphe au-dessus de $S$ qui est $q$-isotrope ou de
 $q$-norme constante \'egale \`a $1$ selon que l'unique sous-groupe \`a un param\`etre contenu dans $B$ est unipotent ou semi-simple.

 En r\'esum\'e, chaque composante connexe de $S$ est {\it lisse} et en chaque point $u  \in S$ il existe un  vecteur $X(u) \in T_{u}M$ tel que toute isom\'etrie locale de $q$ restreinte au sous-ensemble invariant $S$ pr\'eserve $X$.

   Nous d\'emontrons  que le champ de vecteurs $X$ est tangent \`a $S$. Pour tout $u \in S$ les transformations lin\'eaires de $T_{u}M$ qui pr\'eservent
   $X(u)$ pr\'eservent le $s$-jet de $\varphi$ et s'int\`egrent donc en des isom\'etries locales de $\varphi$ qui fixent tous les points de la g\'eod\'esique issue de $u$ en la direction $X(u)$. Comme l'isotropie de l'action de $\cal G$ est triviale sur $M \setminus S$, ceci implique que cette g\'eod\'esique est enti\`erement contenue dans $S$ et donc $X(u) \in T_{u}S$. Nous avons prouv\'e en m\^eme temps que chaque composante connexe de $S$ est de dimension complexe au moins $1$.
      
 On   d\'emontre maintenant que    chaque composante connexe de $S$ est de dimension complexe  $2$ (de codimension $1$ dans $M$). Supposons le contraire et consid\'erons que $S$ est une telle composante connexe de dimension complexe $1$. Associons \`a chaque \'el\'ement de l'alg\`ebre de Lie de champs de Killing $\cal G$  au voisinage d'un point $u \in S$ sa valeur en $u$. Ceci donne un morphisme
 d'espaces vectoriels d\'efini sur $\cal G$ et \`a valeurs dans $T_{u}M$. Comme $S$ est invariant par l'action de $\cal G$, l'image de notre morphisme est incluse dans la droite $T_{u}S$. Ceci implique qu'il existe au moins deux champs de Killing lin\'eairement ind\'ependents dans le noyau et donc 
  le groupe d'isotropie en $u$ serait  n\'ecessairement de dimension au moins $2$. Dans le cas o\`u ce groupe d'isotropie est de dimension $2$, le stabilisateur du $s$-jet de $g$
  en $u$ est de dimension $2$, ce qui montre que ce $s$-jet poss\`ede  une orbite de dimension $2$. Ce cas a \'et\'e exclus par le lemme~\ref{dimension}. Si le groupe d'isotropie est
  de dimension $3$, il est n\'ec\'essairement isomorphe \`a $PSL(2, \CC)$ : il vient que $\cal G$ est isomorphe \`a l'alg\`ebre de Lie de $PSL(2, \CC)$ qui est
  unimodulaire (car semi-simple). Ceci est en contradiction avec la conclusion de la proposition~\ref{unimodulaire}.
  \end{demonstration}
  
  Nous savons \`a pr\'esent que $S$ est de dimension complexe $2$. Dans ce cas une isom\'etrie locale qui est triviale en restriction \`a $S$ admet un $1$-jet trivial
  en chaque point de $S$ et elle est, par cons\'equent, triviale. La restriction de l'alg\`ebre de Lie $\cal G$  \`a $S$ est donc un isomorphisme d'alg\`ebres de Lie.

  Rappelos que le champ de vecteurs $X$ tangent \`a une composante connexe de $S$ construit pr\'ec\'edemment   est $q$-isotrope ou bien de $q$-norme constante \'egale \`a $1$ selon que le sous-groupe
 \`a un param\`etre contenu dans $B$ qui le fixe est respectivement unipotent ou semi-simple. Cette consid\'eration s\'epare les deux cas qu'on \'etudie dans la suite.

  \subsection{Isotropie unipotente}

  Pla\c cons-nous d'abord  dans le cas o\`u {\it il existe dans  l'image de $\varphi^{(s)}$  au moins une orbite $O_{1}$ de dimension $2$, contenue dans $\bar O \setminus O$ et dont le stabilisateur contient un sous-groupe
  \`a un param\`etre  unipotent.}  
  
  Nous \'etudions   une  composante connexe de $S$ (que l'on note encore $S$) o\`u le $s$-jet de $\varphi$ est dans $O_{1}$. Il existe alors sur un rev\^etement  fini de $S$ un champ de vecteurs $q$-isotrope $X$ pr\'eserv\'e par la restriction \`a $S$ de tout champ de Killing local.   Le raisonnement suivant \'etant local on consid\`ere que $X$ est un champ
  de vecteurs d\'efini directement sur $S$. 
  
    Dans ce cas nous allons d\'emontrer la
    
    \begin{proposition}   \label{Kodaira}
    
    La surface compacte connexe lisse $S$ est totalement g\'eod\'esique et d\'eg\'en\'er\'ee pour la m\'etrique riemannienne holomorphe $q$. Le feuilletage holomorphe
    engendr\'e par le noyau de la restriction de $q$ \`a $S$  est transversalement riemannien.   Par cons\'equent, $S$ est ou bien un tore complexe,  ou bien une surface de
    Kodaira primaire (fibr\'e principal en courbes elliptiques sur une courbe elliptique). \end{proposition}
    
    \begin{demonstration}
    
   La r\'estriction de la m\'etrique riemanienne holomorphe \`a $S$ est n\'ecessairement d\'eg\'en\'er\'ee.
    En effet, supposons pour un instant le contraire et d\'esignons par $Z$ l'unique  champ de vecteurs tangent \`a $S$,
   colin\'eaire \`a la deuxi\`eme direction isotrope de la restriction de $q$ \`a l'espace tangent
   \`a $S$ (la droite engendr\'ee par $X$ \'etant l'autre direction isotrope) et tel que $q(X,Z)=1$ et par $H$ l'unique champ de vecteurs de norme 
   $1$ orthogonal au plan engendr\'e par $X$ et $Z$ et tel que $vol(H,X,Z)=1$. Toute isom\'etrie locale qui relie deux points
   de $S$ pr\'eserve n\'ecessairement $X$ et donc $Z$ et $H$. 
   Nous avons expliqu\'e un peu plus haut que le point  $u$ admet un voisinage ouvert $V$ dans $M$ tel que toute isom\'etrie de $\varphi$
    proche de l'identit\'e d\'efinie dans un (petit) ouvert connexe de $V$  se prolonge \`a tout l'ouvert $V$. Consid\'erons une telle isom\'etrie qui relie deux points de $V  \setminus S$ (une telle
    isom\'etrie existe car $\varphi$ est localement homog\`ene sur $M \setminus S$) et qui se prolonge donc sur tout
    l'ouvert $V$. Cette  isom\'etrie respecte n\'ecessairement $S$ car le pseudo-groupe des  isom\'etries locales pr\'eserve l'ouvert $M \setminus S$. Toutes les isom\'etries locales qui relient deux points de $S$ pr\'eservent le champ de vecteurs transverse $H$ et donc \'egalement le temps $t$
    du flot g\'eod\'esique de $H$. On vient de voir que  dans l'ouvert 
    $V$ l'action de $\cal G$  fixe chaque feuille du  feuilletage local de dimension $2$ donn\'e par l'image de $S$ par le flot 
    g\'eod\'esique de  $H$ (chaque feuille $exp_{S}(tH)$, \`a $t$ fix\'e, est stabilis\'ee). L'action de $\cal G$ n'est donc pas transitive sur  
     $M \setminus S$, ce qui est absurde.

   Il reste que  la restriction de $q$ \`a l'espace tangent \`a $S
   $ est d\'eg\'en\'er\'ee et  cet espace tangent n'est rien d'autre que l'orthogonal $X^{\bot}$ du champ de vecteurs $X$.  Nous montrons que la surface $q$-d\'eg\'en\'er\'ee $S$
   est {\it totalement g\'eod\'esique}
   
   Le champ de vecteurs $X$   \'etant  d\'efini seulement  sur  $S$, on prendra soin de  consid\'erer sa
   d\'eriv\'ee covariante uniquement  le long de  champs de vecteurs $W$ tangents \`a $S$ (autrement dit, contenus dans $X^{\bot}$). Comme la $q$-norme de $X$ est constante, nous avons  d\'ej\`a que pour tout $W \in X^{\bot}$~: $2 \cdot q(\nabla_{W}X,X)= W \cdot q(X,X)=0.$ Le champ de plans $X^{\bot}$ est donc stable par l'op\'erateur $\nabla_{\cdot}X$, que l'on peut interpr\'eter comme une section
   au-dessus de $S$ du fibr\'e $End(X^{\bot})=End(TS)$     des endomorphismes de $X^{\bot}$ .  
   
   Constatons d'abord  que le champ $X$ est g\'eod\'esique.  Pour s'assurer que $X$ est bien g\'eod\'esique consid\'erons  l'action du champ de Killing (unipotent) qui fixe un point $u$ de $S$. Si $H(u) \in X^{\bot}(u)$
   est un vecteur de $q$-norme unitaire, la diff\'erentielle du temps $t$ du flot de champ de Killing envoie le vecteur $H(u)$ sur $H(u) +t \cdot X(u)$. Comme ce flot
   pr\'eserve la connexion $\nabla$ et le champ de vecteurs $X$ il vient que $\nabla_{H(u)}X = \nabla_{H(u) +t \cdot X(u)}X$, ce qui implique que $\nabla_{X}X$ s'annule
   au point $u$.

   Par cons\'equent l'op\'erateur $\nabla_{\cdot}X$ contient le champ $X$ dans son noyau.
   
   L'autre valeur propre
   de l'op\'erateur $\nabla_{\cdot}X$ (qui est constante sur $S$ car le pseudo-groupe des isom\'etries locales de $q$ qui pr\'eservent $X$ agit transitivement sur S)
   est n\'ecessairement  nulle : dans le cas contraire l'op\'erateur $\nabla_{\cdot}X$ serait diagonalisable et tout champ de Killing devrait pr\'eserver la d\'ecomposition de
   $X^{\bot}$ en espaces propres, or ceci n'est pas r\'ealis\'e pour notre champ de Killing unipotent dont la diff\'erentielle ne fixe aucune autre droite de $X^{\bot}$ \`a part celle engendr\'ee par $X$. Il reste que le champ 
   d'endomorphismes $\nabla_{\cdot}X$ est nilpotent (d'ordre au plus $2$)~: l'image de $  \nabla_{\cdot}X$ est incluse dans le noyau de   $  \nabla_{\cdot}X$. 
     Deux cas se pr\'esentent : ou bien  $\nabla_{\cdot}X$  est nul,  ou bien    le noyau de $\nabla_{\cdot}X$ et l'image de $\nabla_{\cdot}X$ co\"{\i}ncident avec la
     droite engendr\'ee par $X$ (l'unique droite de $X^{\bot}$ invariante par toutes les isom\'etries locales). Dans les deux cas le calcul suivant est valide pour tous les
     champs de vecteurs    locaux $W_{1}$ et $W_{2}$ tangents \`a $S$~:
     $q(\nabla_{W_{1}}   W_{2}, X)= W_{1} \cdot q(W_{2}, X) - q(\nabla_{W_{1}}  X, W_{2})=0,$ le deuxi\`eme terme du membre de droite de l'\'egalit\'e \'etant nul car
    $ \nabla_{W_{1}}  X   $ est contenu dans l'image de $\nabla_{\cdot}X$ et donc colin\'eaire \`a $X$ (tandis que $W_{2} \in X^{\bot}$). Ceci montre que     $  \nabla_{W_{1}}   W_{2}     \in X^{\bot}$, et que
    $S$ est totalement g\'eod\'esique.
     
     Comme $q$ est d\'eg\'en\'er\'ee en restriction \`a la surface totalement g\'eod\'esique $S$, le feuilletage engendr\'e par le champ $q$-isotrope
     $X$ est transversalement riemannien~\cite{Ze}. 
     
     Le feuilletage engendr\'e par le champ de vecteurs non singulier $X$ admet une structure transverse model\'ee sur le groupe de
     Lie $\CC$. On en d\'eduit facilement (voir~\cite{Mo} pour la th\'eorie g\'en\'erale et \cite{Bos} pour le cas des surfaces complexes) que le
     rev\^etement universel de $S$ est biholomorphe \`a $\CC^2$ et qu'il existe seulement deux cas possibles pour $S$ : ou bien, $S$
     est un tore complexe, quotient de $\CC^2$ par un r\'eseau de translations et $X$ est un champ de vecteurs constant, ou bien $S$ est une surface de Kodaira primaire, fibr\'e principal en courbes elliptiques sur une courbe elliptique et le champ $X$  engendre la fibration principale.
 \end{demonstration}
        
        Analysons  maintenant l'action de l'alg\`ebre de Lie $\cal G$ des champs de Killing de $\varphi$ au voisinage d'un point $u$ de $S$. La restriction \`a $S$ de chaque  \'el\'ement de $\cal G$ donne un champ de vecteurs (tangent \`a $S$) d\'efini au voisinage de $u$ dans $S$
        dont le flot  pr\'eserve $X$ et donc, en particulier, le feuilletage (transversalement riemannien) $\cal F$ d\'efini par $X$. D\'esignons par $\cal H$ l'id\'eal de $\cal G$   form\'e
        par les \'el\'ements de $\cal G$ dont la restriction \`a $S$ agit trivialement sur la transversale de $\cal F$. 
        Les \'el\'ements de $\cal H$ fixent chaque feuille de  $\cal F$ et ils commutent avec $X$ : ils sont de la forme $f \cdot X$ avec $f$ fonction holomorphe constante sur les orbites de $X$ (en particulier, les \'el\'ements de $\cal H$ sont $q$-isotropes sur $S$).  Remarquons que l'alg\`ebre de Lie $\cal H$ est de dimension complexe 
        $2$ : en effet,  $\cal H$ ne peut \^etre de dimension $3$ ( sinon le groupe d'isotropie de $\varphi$ en $u$ serait de dimension au moins $2$) et le quotient $\cal G$/$\cal H$ qui agit non trivialement sur la transversale de $\cal F$ est n\'ecessairement de dimension $1$ isomorphe \`a l'alg\`ebre de Lie $\CC$ (agissant par translation).  Comme $\cal F$ est transversalement riemannien, le flot de $X$ (comme le flot de tout
        champ de  vecteurs tangent au feuilletage) pr\'eserve la restriction de $q$ \`a $S$.

        On peut  choisir sur un voisinage ouvert $U$ de $u$ dans $S$ un champ de vecteurs holomorphe $H$ de $q$-norme constante \'egale \`a $1$ et tel que
        $\lbrack X, H \rbrack=0$ (il suffit de d\'efinir $H$ de $q$-norme constante \'egale \`a $1$ sur une petite transversale \`a $\cal F$ et de le transporter par le flot
        du champ  $X$ qui pr\'eserve la restriction de $q$ \`a $S$). D\'efinissons sur  un voisinage de $u$ dans $S$ un syst\`eme de  coordonn\'ees $(x,h)$ centr\'e en $u$  et tel que  $\frac{\partial}{\partial x}=X$ et $\frac{\partial}{\partial h}=H$.   Dans ces coordonn\'ees l'expression locale de la  forme quadratique $q$ restreinte \`a $S$ est $dh^2$ et
       les \'el\'ements de $\cal H$ restreint \`a $S$ sont de la forme $f(h)  \frac{\partial}{\partial x}$ (car ils pr\'eservent $\frac{\partial}{\partial x}$
       et $dh^2$).
       Par ailleurs l'alg\`ebre de Lie        $\cal G$/$\cal H$       est engendr\'ee par un \'el\'ement qui s'exprime $\frac{\partial}{\partial h} + l(h) \frac{\partial}{\partial x}$, avec $l$ une fonction
       holomorphe d\'efinie au voisinage de $0$ dans $\CC$. Nous avons que 
       $\lbrack \frac{\partial }{\partial h} + l(h) \frac{\partial}{\partial x}, f(h) \frac{\partial}{\partial x} \rbrack  = f'(h) \frac{\partial}{\partial x}.$

       Pour comprendre la structure de l'alg\`ebre de Lie $\cal G$  qui agit par isom\'etries affines  pour la restriction de la connexion $\nabla$
       \`a $S$, nous allons pr\'eciser la structure locale de cette connexion sur $S$.

       \begin{proposition}    \label{courbure}
       
     i).   Si $R$ est le tenseur de courbure de $(S, \nabla )$ alors $R(X,H)X=0$ et $R(X,H)H=\gamma X$, o\`u $\gamma$ est un nombre complexe.
     
     ii). La restriction de $\nabla$ \`a $S$ est localement sym\'etrique. De plus $(S, \nabla)$ est localement isom\'etrique ou bien \`a la connexion
     canonique  du groupe affine de la droite complexe si $\gamma \neq 0$, ou bien \`a la connexion
     canonique de $\CC^2$ si $\gamma=0$.
     
     \end{proposition}

     Avant de passer \`a la preuve rappelons que  le rev\^etement universel $AG$ du groupe affine de la droite complexe est un groupe de Lie complexe de dimension $2$
     qui peut \^etre vu comme l'ensemble des couples $(a,b) \in \CC^2$ muni de la multiplication $(a,b) \cdot (a',b')= (a+a', exp(a)b'+b).$ Ce groupe admet (comme tout groupe de Lie) une unique  connexion lin\'eaire holomorphe, bi-invariante, sans torsion, complete et localement sym\'etrique. Le groupe d'isom\'etries de cette connexion est form\'e par
     les translations \`a droite et \`a gauche et il est, par cons\'equent, isomorphe au produit $AG \times AG$~\cite{Zeg}.
     
     Avant de passer \`a la preuve rappelons que si $h$ est l'\'el\'ement de l'alg\`ebre de Lie du groupe affine qui engendre les homoth\'eties (le temps $T$ de
     du  flot de $h$ agissant donc sur la droite complexe comme $z \to exp(T)z$) et $x$ est  l'\'el\'ement de l'alg\`ebre de Lie qui engendre les translations (le temps $T$ du flot 
     \'etant $z \to z + T$) nous avons $\lbrack h, x \rbrack= -x$. 
     
   La connexion canonique d'un groupe de Lie est d\'efinie, en g\'en\'eral, au niveau de l'alg\`ebre de Lie par la relation $\nabla_{u}v= \frac{1}{2} \lbrack u, v \rbrack$; ce qui donne
     pour la courbure le tenseur $R(u,v)v= \frac{1}{4} \lbrack v, \lbrack u, v \rbrack \rbrack$. 
     
     Dans le cas particulier du groupe affine  il vient que $R(x,h)x=0$ et $R(x,h)h= - \frac{1}{4}x$.
     Si l'on pose $X=x$ et $H=h$, ceci ressemble formellement, du moins dans le cas $\gamma = -a^2= - \frac{1}{4}$, aux  relations figurant au  point i) de la proposition \ref{courbure}. 
     
     Les arguments suivants inspir\'es de~\cite{Zeg} (partie 8) montrent que la valeur du param\`etre  $a$ n'est pas relevante pour la g\'eom\'etrie de la connexion,
      tant que $a$ reste non nul.

     \begin{demonstration}

     Remarquons d'abord que $\nabla_{H}X$ ne d\'epend pas du champ de vecteurs $H$ de $q$-norme unitaire chosi. En effet, si $H'$ est
     un autre champ de vecteurs local tangent \`a $S$ et de $q$-norme constante \'egale \`a $1$, alors $H'=H +f \cdot X$, pour une certaine
     fonction holomorphe locale $f$ et comme $X$ est g\'eod\'esique, $\nabla_{H} X= \nabla_{H'}X$. Le champ de vecteurs $\nabla_{H}X$ est donc invariant
     par l'action de $\cal G$ et, comme  cette action est transitive sur $S$, ceci  implique qu'il existe un nombre complexe $a$ tel que $\nabla_{H}X=aX$.
     
     Comme $\lbrack X, H \rbrack=0$, nous avons que $R(X,H)X= \nabla_{X}  \nabla_{H}  X - \nabla_{H} \nabla_{X}X=\nabla_{X}(aX)=0$.
     
     Pour la deuxi\`eme formule sur la courbure remarquons que le terme $R(X,H)H$ ne d\'epend pas du choix de $H$ car si $H'=H +f \cdot X$,
     alors la premi\`ere \'egalit\'e implique que $R(X,H)f \cdot X=0$ et donc $R(X,H)H=R(X,H')H'$. Comme pr\'ec\'edemment le champ de vecteurs $R(X,H)H$ est alors
     invariant par l'action de $\cal G$,      ce qui implique que $R(X,H)H= \gamma X$ pour un certain nombre complexe $\gamma$. Un calcul direct montre  que $\gamma =-a^2$.

      La deuxi\`eme assertion de la proposition est  prouv\'ee dans~\cite{Zeg} sous l'hypoth\`ese suppl\'ementaire $\nabla_{H}H=0$ qui est utilis\'ee
      pour montrer que la courbure est parall\`ele (connexion localement sym\'etrique sur $S$). Nous montrons que
      dans notre situation  nous pouvons nous passer de l'hypoth\`ese faite sur $H$. Remarquons d'abord que $H$ \'etant de $q$-norme constante, $\nabla_{H}H$ est orthogonal \`a
      $H$. Comme $\nabla_{H}H \in X^{\bot}$ ceci implique qu'il existe une fonction holomorphe $g$ telle que $\nabla_{H}H =g \cdot X$.
      
      Rappelons que la d\'eriv\'ee du tenseur $R$ est donn\'ee par la formule : $$\nabla R(A,B,C,D)=\nabla_{A}(R(B,C)D) -R(\nabla_{A}B,C)D -R(B, \nabla_{A}C)D-R(B,C) \nabla_{A}D,$$ o\`u les   vecteurs $A,B,C, D$ prennent les valeurs $X$ ou $H$. Ce n'est que dans le cas o\`u trois des vecteurs A,B,C,D valent $H$ que l'hypoth\`ese suppl\'ementaire est
      utilis\'ee dans~\cite{Zeg}. Nous traitons ici ce  cas et pour le reste de la preuve nous renvoyons \`a~\cite{Zeg} (proposition 8.4 et proposition 9.2).

1)   $ \nabla R(H,X,H,H)=$ $$ \nabla_{H}R(X,H)H-R(\nabla_{H}X,H)H- R(X, \nabla_{H}H)H- R(X,H) \nabla_{H}X=$$ $$\nabla_{H} \gamma X -R(aX, H)H -R(X,gX)H-R(X,H)aX=
      a \gamma X -a \gamma X =0.$$

2)     $ \nabla R(H,X,H,H)=$ $$ \nabla _{H}R(X,H)H -R(\nabla_{H}X,H)H -R(X, \nabla_{H}H)H -R(X,H)\nabla_{H}H =$$ $$ \nabla_{H} \gamma X - R(aX,H)H - R(X,gX)H - R(X,H)gX=
     a \gamma X- a \gamma X -0 -0 =0.$$
     
 3)    $\nabla R(H,H,H,X)=$ $$ \nabla_{H} R(H,H)X - R( \nabla_{H}H,H)X - R(H, \nabla_{H}H)X-R(H,H) \nabla_{H}X=$$
     $$0-R(gX,H)X -R(H,gX)X-0=0.$$
      \end{demonstration}

       Revenons \`a pr\'esent \`a
      l'alg\`ebre de Lie $\cal G$. Dans le cas o\`u $\gamma=0$ notre alg\`ebre de Lie se plonge dans l'alg\`ebre de Lie du groupe des transformations affines de $\CC^2$ qui pr\'eservent
      la forme quadratique $dh^2$ et le  champ de vecteurs $X$ qui est isotrope, parall\`ele, et non trivial. Cette alg\`ebre  est engendr\'ee par $\frac{\partial}{\partial x}, \frac{\partial}{\partial h}, h\frac{\partial}{\partial x}$.
        Il s'agit de l'alg\`ebre de Lie du groupe Heisenberg car $\frac{\partial}{\partial x}$ est dans le centre et $\frac{\partial}{\partial x}= \lbrack \frac{\partial}{\partial h}, h \frac{\partial}{\partial x} \rbrack$. L'alg\`ebre de Heisenberg  est unimodulaire (car nilpotente) et cette situation a d\'ej\`a \'et\'e analys\'ee et \'elimin\'ee comme contradictoire.      
        
        Supposons maintenant que $\gamma \neq 0$. Dans ce cas $\cal     G$ se plonge dans l'alg\`ebre de Lie du produit $AG \times AG$  (comme
        $AG$ est simplement connexe et complet, toute isom\'etrie locale proche de l'identit\'e de $\nabla$ se prolonge en une isom\'etrie globale~\cite{Amo}, \cite{Gro}, \cite{No}). 
        
        Nous allons montrer que l'alg\`ebre de Lie  $\cal G$ admet un centre non trivial. Pour cela
	nous prouvons d'abord que $X$ est un champ de Killing pour la restriction de $\nabla$
	\`a $S$.
        
       Rappelons que, d'apr\`es la proposition \ref{Kodaira}, la surface   $S$ est biholomorphe \`a un tore complexe
       ou bien \`a une surface de Kodaira primaire.
      
      Consid\'erons d'abord le cas o\`u  $S$ est un tore complexe. Il vient que le fibr\'e holomorphe tangent \`a $S$ est holomorphiquement trivial et que le champ de vecteurs local $H$ de norme constante \'egale \`a $1$ peut \^etre choisi comme \'etant globalement d\'efini sur  $S$.
   De m\^eme le champ de vecteurs $\nabla_{H}H$ sera globalement d\'efini et comme $H$ est de $q$-norme constante, $\nabla_{H}H$ est orthogonal \`a $H$, donc
   en tout point collin\'eaire \`a $X$. Comme le fibr\'e en droites d\'efini par $X$ est trivial (car $X$ non singulier), on a que $\nabla_{H}H=\lambda X$, pour
   une certaine constante $\lambda \in \CC^*$.
   
   Comme le flot de $X$ pr\'eserve $H$, les relations $\nabla_{X}X=0$, $\nabla_{H}X= \nabla_{X}H=a X$ et $\nabla_{H}H=\lambda X$ impliquent
   que le flot de $X$ pr\'eserve la connexion $\nabla$ de $S$.  
   
   Une mani\`ere plus directe de conclure est de dire que toute connexion lin\'eaire holomorphe
   sur un tore complexe est invariante par les translations (pour une description des connexions
   affines sur les tores complexes le lecteur pourra consulter~\cite{IKO}).
   
   Dans le cas o\`u $S$ est une surface de Kodaira primaire, nous utilisons la classification
   des connexions lin\'eaires  holomorphes sur ce type de surface faite par A. Vitter dans~\cite{Vi} (page 239). Cette description montre, en particulier, que sur une surface de Kodaira primaire toute connexion lin\'eaire  holomorphe est invariante par la fibration principale. Dans notre cas, il vient que $X$ est bien un champ de Killing pour la restriction de $\nabla$ \`a $S$.
   
   Dans le mod\`ele local constitu\'e par le groupe affine, l'\'el\'ement $X$ correspond alors \`a un \'el\'ement $x$ de l'alg\`ebre de Lie du produit $AG \times AG$. Comme
   les \'el\'ements de $\cal G$ pr\'eservent $X$, l'alg\`ebre de Lie $\cal G$ est une sous-alg\`ebre de Lie de dimension $3$ du commutateur de $x$ dans l'alg\`ebre de Lie 
   du produit $AG \times AG$. Comme, par ailleurs,  le commutateur de chaque \'el\'ement $x$ est une alg\`ebre de Lie de dimension au plus $3$, il vient que $\cal G$
   co\"{\i}ncide avec le commutateur de $x$ et poss\`ede donc $x$ comme \'el\'ement central non trivial. 
   
   On vient de prouver que $X$ est la restriction \`a $S$ d'un champ de Killing local qui se trouve dans le centre de  $\cal G$. Comme ce champ de Killing est invariant par l'action de $G$ sur lui-m\^eme, il fournit un champ de Killing globalement
  d\'efini sur $M \setminus S$ et de norme constante (n\'ecessairement \'egale \`a $0$ car ce champ de Killing se prolonge
  sur un voisinage ouvert de $S$ en le champ $q$-isotrope $X$).
  
  Analysons la situation au voisinage d'un point $u$ de $S$. Pour cela, consid\'erons un voisinage $U$ de $u$ dans $M$ qui satisfait
  la propri\'et\'e de prolongement de champs de Killing et consid\'erons dans $U \setminus S$ un champ de Killing $K_{1}$ qui correspond
  \`a l'\'el\'ement central de $\cal G$. Ce champ de Killing est de norme constante sur $U \setminus S$ car invariant par l'action transitive de $\cal G$  et se prolonge sur $U \cap S$ en un champ colin\'eaire \`a $X$. Nous avons donc que $K_{1}$ est bien d\'efini sur $U$ et
  $q$-isotrope. Rappelons  que la section $Y$ de $S^n(TM)$ est \'egalement  invariante par l'action de $\cal G$ 
  et (par construction) $q$-isotrope. Nous avons donc que  $q(K_{1}^{\otimes n}, Y )$ est une fonction constante, o\`u l'on a not\'e encore par $q$
  la forme quadratique induite  sur $TM^{\otimes n}$ par la m\'etrique riemannienne holomorphe. Comme $Y$ s'annule sur $S$ cette
  fonction est nulle. Nous pouvons donc conclure que si l'on choisit en un point $v \in U \setminus S$ le vecteur $\tilde Y(v)$ tel que
  $\tilde Y(v)^{\otimes n}=Y(v)$ alors les vecteurs $\tilde Y(v)$ et $ K_{1}(v)$       sont $q$-isotropes et $q$-orthogonaux, ce qui implique
  qu'ils sont colin\'eaires : il existe donc une fonction m\'eromorphe $g$, d\'efinie sur $U$, telle que  que $K_{1}^{\otimes n} =g \cdot Y$. L'invariance de $K_{1}^{\otimes n}$    et de $Y$ par l'action transitive de $\cal G$  implique que $g$ est constante et 
  qu'il existe donc $\mu \in \CC^{*}$ tel que $K_{1}^{\otimes n} =\mu Y$ sur $U$. La contradiction recherch\'ee vient du fait que $Y$ s'annule sur $S$, tandis que
  le champ de Killing $K_{1}$ est non identiquement nul sur $U \cap S$ (car nous avons vu que  la restriction \`a $S$ est un isomorphisme). En fait, en tant que champ de Killing
  $q$-isotrope en dimension $3$, le champ de vecteurs  $K_{1}$ est m\^eme non singulier~\cite{D2} (lemme 6.7).
        
      \subsection{Isotropie   semi-simple}
     
    Il reste \`a r\'egler les cas o\`u l'isotropie est semi-simple et de dimension $1$. 
  Pla\c cons-nous    dans le dernier cas \`a r\'egler  o\`u {\it en dehors de $O$, il n'existe dans  l'image de $\varphi^{(s)}$   que des orbites de dimension $2$ et dont le stabilisateur contient un sous-groupe
  \`a un param\`etre semi-simple.}

  Nous \'etudions   une  composante connexe de $S$ (que l'on note encore $S$) o\`u le $s$-jet de $\varphi$ est dans une orbite $O_{1}$ de dimension $2$ et dont le stabilisateur contient un sous-groupe \`a un param\`etre semi-simple. Ce stabilisateur fournit en chaque point $u$ de $S$ un champ de Killing s'annulant en $u$ et dont la diff\'erentielle
  fixe un unique vecteur $q$-unitaire de $T_{u}S$.  Il existe alors sur $S$ un champ de vecteurs tangent $X$ de $q$-norme  constante \'egale \`a $1$  pr\'eserv\'e par la restriction \`a $S$ de tout \'el\'ement de $\cal G$. 
  
    La m\'etrique riemannienne holomorphe $q$ est n\'ecessairement  d\'eg\'en\'er\'ee sur $S$
     car sinon il existe sur $S$ un deuxi\`eme champ de vecteurs $H$  d\'efini par les relations $q(X,H)=0$ et $q(H)=1$. Le champ $H$ est bien d\'efini seulement \`a signe pr\`es, car sur la droite orthogonale \`a $X$ il existe deux vecteurs (oppos\'es) $q$-unitaires. Faisons un choix local et
     optons pour l'un de ces vecteurs au voisinage d'un point $s$ de $S$ : ceci permet de d\'efinir 
     un champ $H$ pr\'eserv\'e par l'action de $\cal G$.

     L'orthogonal (par rapport \`a $q$) de l'espace tangent \`a $S$ est engendr\'e par un  unique champ de vecteurs $Z$ de norme $1$ tel que $vol(X,H,Z)=1$. Les champs $H$ et $Z$ sont donc \'egalement pr\'eserv\'es  par l'action au voisinage du point $s$ de la restriction de
     $\cal G$ \`a $S$. On conclut alors qu'au voisinage d'un point de $S$ l'action de $\cal G$  fixe chaque feuille  du feuilletage obtenu  en poussant $S$
     par le flot g\'eod\'esique de $Z$  : ceci implique que $\cal G$  n'agit  pas transitivement sur $M \setminus S$, ce qui est absurde.
     
     Il reste \`a analyser le cas o\`u $q$ est d\'eg\'en\'er\'ee en restriction \`a $S$. Le noyau de $q$ d\'efinit alors en restriction \`a $S$ un feuilletage en droites non singulier
     transverse au champ de vecteurs $X$.  Nous sommes \`a nouveau en pr\'esence d'un feuilletage holomorphe transversalement riemannien. En effet, d\'esignons par  $\omega$ la $1$-forme diff\'erentielle holomorphe sur $S$ qui vaut $1$ sur le champ de vecteurs $X$ et  dont le noyau co\"{\i}ncide avec
     celui de $q$. Il vient  que $\omega^2= q$ et que $\omega$ est ferm\'ee;  nous venons d'utiliser  un  r\'esultat classique  qui affirme  que  sur 
     les  surfaces complexes compactes, les $1$-formes diff\'erentielles holomorphes sont n\'ecessairement ferm\'ees (pour une preuve de ce fait le lecteur
     pourra consulter, par exemple, la r\'ef\'erence ~\cite{Bru}).  Par cons\'equent, le noyau de $\omega$
     d\'efinit bien un feuilletage transversalement riemannien $\cal F$ et $S$ est n\'ecessairement totalement g\'eod\'esique~\cite{Ze}. Il r\'esulte, en utilisant la classification des surfaces
     poss\`edant une connexion lin\'eaire holomorphe~\cite{W} que $S$ ne peut \^etre qu'un tore complexe ou une surface de Hopf. N\'eanmoins le raisonnement  local suivant
     n'utilise pas  la structure analytique globale de $S$.
     
     \begin{proposition}  \label{champ}
     
     Il existe un champ de vecteurs holomorphe $\tilde Y$ tel que $Y = \tilde Y^{\otimes n}$. De plus, le $1$-jet de $\tilde Y$ en un point de $S$ est non trivial.
     \end{proposition}
     
     \begin{demonstration}
     
     Analysons l'action de l'alg\`ebre de Lie des champs de Killing de $\varphi$ au voisinage d'un point $u$ dans $S$.  Dans  un voisinage (suffisamment petit) $U$ de $u$ dans $S$, 
     il existe   un champ de vecteurs holomorphe tangent $H$ qui est $q$-isotrope et non singulier dans $U$. En chaque point de $U$ il existe alors un unique vecteur $Z$ (transverse \`a $U$) uniquement d\'etermin\'e
     par les relations $q(Z)=0, q(X,Z)=0$ et $q(H,Z)=1$. Le flot g\'eod\'esique de $Z$ permet de d\'efinir un feuilletage au voisinage de $U$ dont les feuilles sont $exp_{U}(tZ)$,
     \`a $t$ fix\'e.
          
    Rappelons qu'un $1$-jet de biholomorphisme local qui relie deux points de $S$ se prolonge en une isom\'etrie locale si et seulement si sa diff\'erentielle pr\'eserve $X$.  Le $1$-jet d'application qui envoie $u$ sur un point $u' \in U$ et dont la diff\'erentielle est triviale dans les 
      rep\`eres $(X(u),H(u),Z(u))$  et  $(X(u'), H(u'), Z(u'))$ se prolonge donc de mani\`ere unique en une isom\'etrie locale de $\varphi =(q, Y)$  qui  envoie le temps $t$ de la g\'eod\'esique issue de $u$ dans la direction $Z(u)$ sur le temps $t$ de la g\'eod\'esique issue de $u'$ dans la direction $Z(u')$. Le pseudo-groupe d'isom\'etries locales de $\varphi$  agit donc transitivement sur chaque feuille $exp_{U}(tZ)$ du feuilletage. Plus pr\'ecis\'ement, $u$ \'etant  fix\'e, quelque soit $u' \in U$ il existe une unique isom\'etrie locale $i_{u'}$ (qui d\'epend
      de mani\`ere holomorphe de $u'$)
      qui relie $u$ et $u'$ et qui envoie la  g\'eod\'esique $exp_{u} (tZ(u))$ sur la g\'eod\'esique $exp_{u'} (tZ(u'))$ en respectant le param\'etrage.
           
    Le tenseur $Y$ est
     enti\`erement  d\'etermin\'e si l'on conna\^{\i}t sa restriction \`a la g\'eod\'esique issue de $u$ dans la direction du vecteur $Z(u)$ (en effet, on obtient alors la restriction de
     $Y$ \`a la g\'eod\'esique $exp_{u'} (tZ(u')$ en transportant avec $i_{u'}$). 
     
     Pour pr\'eciser ceci consid\'erons sur un voisinage ouvert de $u$ dans $U$, un syst\`eme de  coordonn\'ees locales $w$, centr\'e en $u$. Un voisinage ouvert de $u$ dans
     $M$ sera alors param\`etre \`a l'aide du flot g\'eod\'esique de $Z$ par les coordonn\'ees $(w,t)$. D\'esignons par $Y'$ la restriction de $Y$ \`a la g\'eod\'esique issue de $u$ dans la direction $Z$. Il vient que  $Y'(t)=Y(0,t) =f(t)(\bar Y(t))^{\otimes n}$, o\`u $f$ est une fonction holomorphe \`a une variable d\'efinie  sur un voisinage de $0$ et qui admet un z\'ero isol\'e
      en $0$ et $\bar Y(t)$ est un vecteur non nul  et $q$-isotrope de $T_{(0,t)}M$ (cette \'ecriture n'est pas unique).  Si $i_{w}$ est l'unique isom\'etrie locale qui relie $0$ \`a $w$ et qui envoie
     la  g\'eod\'esique $exp_{0} (tZ(0))$ sur la g\'eod\'esique $exp_{w} (tZ(w))$ en respectant le param\'etrage, alors $Y(w,t)=f(t)(di_{w} \cdot \bar Y(t))^{\otimes n}.$

    Si  $\mu t^l$, avec $\mu$ nombre complexe non nul  et $l$ entier strictement positif, est le premier jet non nul  de $f$ en $0$, alors le premier jet non nul de $Y'$ \`a l'origine (qui s'interpr\`ete comme un polyn\^ome homog\`ene d\'efini sur la droite $Z(u)$ et \`a valeurs dans l'espace vectoriel $T_{u}M^{\otimes n}$) est de la forme $\mu t^l \bar Y(0)^{\otimes n}= \mu t^l (aX(0)+bH(0) +cZ(0))^{\otimes n}$, avec $a,b,c$ des nombres complexes dont au moins un est non nul. Nous d\'emontrons que l'invariance de cette expression par le flot du champ de Killing semi-simple $K$ qui fixe $u$
    (flot qui stabilise  la g\'eod\'esique issue de $u$ dans la direction $Z(u)$ sans pr\'eserver le param\'etrage et fixe $X(u)$) implique que  le premier jet non nul de $f$ \`a l'origine est  de la forme 
    $ \mu t^n$, avec $\mu$ nombre complexe non nul et $a=b=0$.
   
    En effet, l'action de la diff\'erentielle en $u$ du temps $T$ de ce  flot est $$T \cdot (X(0), H(0), Z(0))=(X(0),T^2H(0),T^{-2}Z(0)),$$  tandis que l'action sur le $l$-jet de $f$ est
    $T \cdot \mu t^l=T^{2l} \mu t^l. $
    Le $l$-jet   de $Y'$ se d\'ecompose en une somme dont les termes sont de la forme
 $\mu t^la^pb^qc^rX(0)^p \otimes  H(0) ^q \otimes Z(0)^r,$ o\`u $p,q,r$ sont des entiers positifs dont la  somme vaut $n$. 
 
    Comme l'action lin\'eaire du temps $T$ du flot de $K$ sur $(T_{0} M)^{\otimes n}$ est diagonalisable dans la base 
$X(0)^p \otimes H(0)^q \otimes Z(0)^r,$ et la valeur propre associ\'ee \`a ce vecteur propre est $T^{2(q-r)}$, l'invariance du $l$-jet de $Y'$
l'action du flot de $K$ implique que  tous les termes de la forme $X(0)^p \otimes  H(0)^q \otimes Z(0)^{ r}$ avec
$p+q+r=n$ qui  ont un co\'efficient non 
nul doivent \^etre   associ\'es \`a la m\^eme valeur propre. Comme $\bar Y(0)$ est $q(0)$-isotrope, ceci montre que la seule possibilit\'e est
$a=b=0$ et $l=n$.

    Ceci implique qu'au voisinage du point $u$ il existe un champ de vecteurs holomorphe  $\tilde Y$ tel que $Y=( \tilde Y)^{\otimes n}$. En effet,  on peut prendre  $\tilde Y(w,t)=l(t) t di_{w} \bar Y(t),$ o\`u l(t) est une ''racine n-\`eme'' 
    de $\frac{f(t)}{t^n}$ dans un voisinage de $0$ dans $\CC$.

    Cette propri\'et\'e \'etant vraie au voisinage de chaque point de $S$ nous avons que (eventuellement sur un rev\^etement fini non ramifi\'e de $M$)  le tenseur $Y$ est la puissance $n$-\`eme d'un champ de vecteurs. 
    
    Consid\'erons ${\tilde Y}'$ la restriction du champ $\tilde Y$ \`a la g\'eod\'esique issue de $u$ dans la direction $Z(u)$, prenons le premier jet non trivial de $\tilde Y'$ en $0$ et
   exprimons  l'invariance de ce jet par le flot du champ de Killing $K$. Nous sommes dans la m\^eme situation que pr\'ec\'edemment et dans le cas particulier $n=1$.
   Le calcul du point i) implique alors que le premier $l$-jet non trivial de $\tilde Y'$ \`a l'origine correspond \`a $l=1$. En particulier, le premier jet de $\tilde Y$ en $u$ est non trivial.
   Comme $\cal G$ agit transitivement sur $S$ en pr\'eservant $Y$ ceci est  vrai pour  tous les points de $S$.
    \end{demonstration}
   
   Nous donnons maintenant un renseignement plus pr\'ecis sur le $1$-jet de $\tilde Y$ en un point de $S$ avec le lemme local suivant :

\begin{lemme}   \label{singularite} Consid\'erons  un voisinage ouvert $U$ de l'origine dans ${\bf C}^3$ muni d'une m\'etrique riemannienne holomorphe $q$.

Soit $\tilde Y$ un champ de vecteurs holomorphe non identiquement nul
$q$-isotrope qui s'annule \`a l'origine. Supposons qu'il existe un champ (de Killing)
holomorphe $K$ dans $U$ qui s'annule \`a l'origine, dont le flot local pr\'eserve \`a la fois $q$ et $Y$ et dont la diff\'erentielle \`a l'origine stabilise un vecteur non isotrope  de $T_{0}\CC^3$ .  Alors  le premier jet
de $\tilde Y$ \`a l'origine  est de la forme  $ \nu z \frac{\partial}{\partial z}$, $\nu \in \CC.$
\end{lemme}

Remarque : rappelons encore une fois que la condition d'invariance  par la diff\'e-rentielle \`a l'origine  d'un champ de vecteurs non isotrope est une condition alg\'ebrique qui porte sur la partie lin\'eaire du champ de Killing $K$ \`a l'origine. Il est \'equivalent de dire que cette partie lin\'eaire est semi-simple.

\begin{demonstration}

i) Le flot du champ $K$ fixe un vecteur $\frac{\partial}{\partial x}$ de $q$-norme \'egale \`a $1$ dans $T_{0}\CC^3$.  Consid\'erons une base $(\frac{\partial}{\partial x},\frac{\partial}{\partial y}, \frac{\partial}{\partial z})$ de l'espace tangent \`a l'origine avec la propri\'et\'e que $\frac{\partial}{\partial x}$ est un vecteur de $q$-norme \'egale \`a $1$, tandis que les vecteurs 
$\frac{\partial}{\partial y}, \frac{\partial}{\partial z}$ sont $q$-isotropes, $q$-orthogonaux \`a $\frac{\partial}{\partial x}$ et $q(\frac{\partial}{\partial y}, \frac{\partial}{\partial z})=1. $

    Dans les coordon\'ees exponentielles fix\'ees par le choix  de la base pr\'ec\'edente de l'espace tangent \`a l'origine le champ de Killing $K$ est lin\'eaire 
    et le temps $T$ de son flot est donn\'e  par la transformation lin\'eaire :
    $T \cdot (x,y,z)=(x, T^{-2}y, T^{2}z).$
    
    Un calcul simple donne la forme des champs de vecteurs holomorphes s'annulant \`a l'origine et invariants par ce flot : $\tilde Y=f(x,y,z) \frac{\partial}{\partial x}    +
    g(x,y,z) \frac{\partial}{\partial y}    + h(x,y,z) \frac{\partial}{\partial z} $ est invariant par le flot de $K$ si et seulement si $f=F(x,yz), g=yG(x,yz)$ et $h=zH(x,yz)$,
    o\`u $F,G$ et $H$ sont des fonctions holomorphes \`a deux variables s'annulant \`a l'origine.
    
       On en d\'eduit les seuls $1$-jets possibles pour des champs de vecteurs s'annulant \`a l'origine et invariants par le flot de $K$ : un tel $1$-jet 
       est n\'ecessairement  de la forme $Y^{1}= \lambda x \frac{\partial}{\partial x}     + \mu y \frac{\partial}{\partial y}  + \nu z \frac{\partial}{\partial z} $, avec $\lambda, \mu, \nu$,
       des nombres complexes eventuellement nuls.
       
       Il s'agit maintenant d'exploiter le fait que $\tilde Y$ est $q$-isotrope. Dans les coordonn\'ees exponentielles consid\'er\'ees le $1$-jet \`a l'origine de la m\'etrique riemannienne
       holomorphe $q$ vaut $q_{0}=dx^2 +dydz$.  Ceci permet d'en d\'eduire le $2$-jet \`a l'origine de la fonction $q(Y)$ : il s'agit de $q_{0}(Y^{1})=
       2 \mu \nu yz + {\lambda}^2x^2$.
          Comme la fonction  $q(Y)$ est identiquement nulle, ce $2$-jet doit  \^etre     trivial et donc $Y^1$ est  de la forme $\nu z \frac{\partial}{\partial z}$   (ou bien
          de la forme $\mu y \frac{\partial}{\partial y}$ ce qui revient au m\^eme) .
   \end{demonstration}
   
   Remarque : dans notre situation on applique le lemme pr\'ec\'edent au voisinage d'un point $u$ de $S$ et pour une base $(\frac{\partial}{\partial x},\frac{\partial}{\partial y}, \frac{\partial}{\partial z}) =(X(u), H(u), Z(u))$ de l'espace tangent en $u$. Comme $\tilde Y$ s'annule sur $S$ et $T_{u}S$ est l'espace vectoriel engendr\'e par 
   $(\frac{\partial}{\partial x},\frac{\partial}{\partial y})$, il vient que le $1$-jet de $\tilde Y$ en $u$ est de la forme $\nu z \frac{\partial}{\partial z}$, avec $\nu$ non nul.

   Le lemme pr\'ec\'edent et la proposition~\ref{champ} impliquent que dans notre situation le $1$-jet de $\tilde Y$  est non nul, de la forme $\nu z \frac{\partial}{\partial z}$, avec
   $\nu$ nombre complexe non nul. La trace de l'op\'erateur   $\nabla_{\cdot} \tilde Y$, vu comme section de $End(TM)$,  est une fonction holomorphe constante sur $M$ qui vaut $\nu$ en les points de
   $S$. Cette fonction est donc une constante non nulle $\nu$. Rappelons que la  trace de $\nabla_{\cdot} \tilde Y$ co\"{\i}ncide avec la divergence du champ de vecteurs $\tilde Y$ par rapport \`a la forme
   volume $vol$ de $q$. Autrement dit, $L_{\tilde Y}vol= \nu vol$, o\`u $L_{\tilde Y}vol$ est la d\'eriv\'ee de Lie de la forme volume par rapport au champ $\tilde Y$. Le fait
   que $\nu$ soit non nulle est en contradiction avec le fait bien connu qu'un champ de vecteurs holomorphe sur une vari\'et\'e compacte doit pr\'eserver les formes volumes holomorphes.
   La preuve de cette propri\'et\'e fait l'objet de la proposition suivante qui ach\`eve la d\'emonstration  de ce dernier cas.
   
   \begin{proposition}
   Le flot de $\tilde Y$ pr\'eserve la forme volume de $q$.
   \end{proposition}
   
   \begin{demonstration}
     
       D\'esignons par $\psi^t$ le temps $t$ du flot de $\tilde Y$ et remarquons que le fibr\'e canonique
  de $M$ \'etant trivial, pour chaque temps $t$ il existe une constante $c(t) \in \CC$ tel que
  $(\psi^t)^*vol=c(t)vol$. Comme le volume r\'eel de $M$, donn\'e par l'int\'egrale de la forme r\'eelle $vol \wedge
  \bar{vol}$, doit  \^etre pr\'eserv\'e par l'action du flot, on a que la valeur absolue du nombre
  complexe $c(t)$ vaut $1$ pour tout  $t$. La fonction $t \to c(t)$ est alors une fonction
  enti\`ere \`a valeurs dans le cercle unit\'e et par le th\'eor\`eme de Liouville elle est n\'ecessairement
  constante, \'egale donc \`a $c(0)=1$. Ceci montre que le flot de $\tilde Y$ pr\'eserve la forme
  volume et par cons\'equent la divergence de $\tilde Y$ est nulle.
  \end{demonstration}

   {\small

Version abr\'eg\'ee du titre : M\'etriques riemanniennes holomorphes.

Titre en anglais :  Holomorphic riemannian metrics on compact threefolds are locally homogeneous

Mots-cl\'es : vari\'et\'es complexes - m\'etriques riemanniennes holomorphes - structures rigides-
pseudo-groupe d'isom\'etries locales.

Classification math. : 53B21, 53C56, 53A55.

\newpage

{\footnotesize

\vspace{2cm}

Sorin Dumitrescu

\rule{4cm}{.05mm}

D\'epartement de Math\'ematiques d'Orsay 

\'Equipe de Topologie et Dynamique

Bat. 425

\rule{4cm}{.05mm}

U.M.R.   8628  C.N.R.S.

\rule{4cm}{.05mm}

Univ. Paris-Sud (11)

91405 Orsay Cedex

France

\rule{4cm}{.05mm}

Sorin.Dumitrescu@math.u-psud.fr

 \end{document}